\documentclass{article}
\usepackage{amssymb}
\usepackage{amsfonts,amsmath}
\usepackage{tikz}
\usepackage[utf8]{inputenc}

\begin{document}

\begin{center}
{\Large Some applications of finite BL-algebras}%
\begin{equation*}
\end{equation*}

\bigskip Cristina Flaut$^{\ast }$, Dana Piciu, Bianca Liana Bercea-Straton%
\begin{equation*}
\end{equation*}
\end{center}

\textbf{Abstract. }{\small In this paper we present an encryption/decryption
algorithm which use properties of finite MV-algebras, we proved that there
are no commutative and unitary rings }$R${\small \ such that }$Id\left(
R\right) =L,${\small where }$L${\small \ is a finite BL-algebra which is not
an MV-algebra and we give a method to generate BL-comets. Moreover, we give
a final characterisation of finite BL-algebra and we proved that a finite
BL-algebra is a comet or MV-algebras which are not chains.}%
\begin{equation*}
\end{equation*}


\begin{equation*}
\end{equation*}

\textbf{1. Preliminaries}%
\begin{equation*}
\end{equation*}

It is known that a commutative ring $R$ for which its lattice of ideals is
isomorphic to an MV-algebra is a direct sums of local Artinian chain rings
with units, see [BN; 09]. Starting from this result, we tried to find
similar characterisation in the case of finite BL-algebras which are not
MV-algebras. But the answer which we found was in the negative sense. In the
paper [NL; 03], authors proved that by using BL-comets, any finite
BL-algebra can be represented as a direct product of BL-comets. In this
paper we proved that there is no commutative and unitary rings $R$ such that
its lattice of ideals, $Id\left( R\right) $, if it is finite, can be
organised as a finite BL-algebra which are not MV-algebra. As a corollary of
this result, we give a characterisation of finite BL-algebras, namely: a
finite BL-algebra is a BL-comet or an unordered MV-algebra, that means an
MV-algebra which is not an MV-chain.

The paper is organised in this introductory part and other three sections.
Section 2 is devoted to present an encryption algorithm based of properties
of an MV-algebra. Section 3 presents the main result of this section,
namely: a finite BL-comet can't be organised as the lattice of ideals of a
commutative and unitary ring $R$. Section 4 gives a method to generate
finite BL-algebras and presents the main result of the paper: there is no
commutative and unitary rings $R$ such that their lattices of ideals, $%
Id\left( R\right) $, if are finite, can be organised as a finite BL-algebra,
which is not an MV-algebra, and, at the end, as a consequence of this
result, we give a characterisation of finite BL-algebras. So, we can
conclude that this paper emphasizes developments of the subject and closes a
problem for the study of finite BL-algebras, regarding their representation
as a lattice of ideals of commutative and unitary ring, but open a direction
to study and characterize infinte BL-algebras.

Let $R$ be a commutative unitary ring. The set $Id\left( R\right) $ denotes
the set of all ideals of the ring $R$. For $I,J\in Id\left( R\right) $, the
following sets are also ideals in $R:$

\begin{equation*}
I+J=<I\cup J>=\{i+j,i\in I,j\in J\}\text{, }
\end{equation*}%
\begin{equation*}
I\otimes J=\{\underset{k=1}{\overset{n}{\sum }}i_{k}j_{k},\text{ }i_{k}\in
I,j_{k}\in J\}\text{, }
\end{equation*}%
\begin{equation*}
\left( I:J\right) =\{x\in R,x\cdot J\subseteq I\}\text{,}
\end{equation*}%
\begin{equation*}
Ann\left( I\right) =\left( \mathbf{0}:I\right) \text{, where }\mathbf{0}=<0>,%
\text{ }
\end{equation*}%
\textbf{\ }and are\textbf{\ }called \textit{sum}, \textit{product}, \textit{%
quotient} and \textit{annihilator} of the ideal $I$.\medskip

\textbf{Remark 1}. ([AF; 92],[AM; 69], [FK; 12])

1)Each nonzero element in a finite commutative unitary ring\ $R$\ is a unit
or a zero divisor.

2) In an Artinian ring every prime ideal is maximal.

3) An Artinian ring is a finite direct product of Artinian local rings.

4) In a commutative ring $R$, the set of non-unit elements is an ideal if
and only if the ring $R$ is local. That ideal is the maximal ideal.\medskip

\textbf{Remark 2.} ([AF; 92],[AM; 69], [FK; 12])

1)Let $R$ be an Artinian commutative ring. Then, each prime ideal is a
maximal ideal.

2) An integral domain $A$ is an Artinian ring if and only if $A$ is a field.

3) An Artinian ring is a finite direct product of Artinian local rings.

4) Let $R$ be a commutative unitary ring.

i) An ideal $M$ of the ring $R$ is maximal if it is maximal, with respect of
the set inclusion, amongst all proper ideals of the ring $R$. From here, it
results that there are no other ideals different from $R$ contained $M$. An
ideal $J$ of the ring $R$ is considered a minimal ideal if it is a nonzero
ideal which contains no other nonzero ideals.

ii) A commutative unitary ring $R$ with a unique maximal ideal is called a
local ring.

iii) We consider $P$ be an ideal in the ring $R,P\neq R$. For $a,b\in R$
such that $ab\in P$, if we have \thinspace $a\in P$ or $b\in P$, therefore $%
P $ is called a prime ideal of $R$.\medskip

\textbf{Definition 3}. ([WD; 39]) \ A \textit{(commutative) residuated
lattice }\ is an algebra $(L,\wedge ,\vee ,\odot ,\rightarrow ,0,1)$ such
that:

(i) $\ (L,\wedge ,\vee ,0,1)$ is a bounded lattice;

(ii) $\ \ (L,\odot ,1)$ is a commutative ordered monoid;

(iii) $\ z\leq x\rightarrow y$\ iff $x\odot z\leq y,$\ for all $x,y,z\in L.$

Property (iii) is called\emph{\ }\textit{residuation}, where $\leq $ is the
partial order of the lattice $\ (L,\wedge ,\vee ,0,1).$

In a residuated lattice we define the following additional operation: for $%
x\in L,$ we denote $x^{\ast }=x\rightarrow 0.$

If we preserve these notations, for a commutative and unitary ring we have
that 
\begin{equation*}
(Id(R),\cap ,+,\otimes \rightarrow ,0=\{0\},1=R)
\end{equation*}%
is a residuated lattice in which the order relation is $\subseteq $, $%
I\rightarrow J=(J:I)$ and $I\odot J=I\otimes J,$ for every $I,J\in Id(R)$,
see [TT; 22]

In a residuated lattice $(L,\wedge ,\vee ,\odot ,\rightarrow ,0,1)$ we
consider the identities:

\begin{equation*}
(prel)\qquad (x\rightarrow y)\vee (y\rightarrow x)=1\qquad \text{ (\textit{%
prelinearity)};}
\end{equation*}

\begin{equation*}
(div)\qquad x\odot (x\rightarrow y)=x\wedge y\qquad \text{ (\textit{%
divisibility)}.}
\end{equation*}

In this paper, by unordered MV-algebra we understand an MV-algebra that is
not chain. By a chain ring $R$ we understand a commutative unitary ring
sucht that its lattice of ideals, $Id\left( R\right) $, is totally ordered
by inclusion.\medskip 

\textbf{Definition 4. }([T; 99])

1) A residuated lattice $L$ is called \textit{a BL-algebra\ }if in $L$ are
verified conditions $(prel)$ and $(div)$.

2) A \textit{BL-chain}\emph{\ }is a totally ordered BL-algebra, that means
it is a BL-algebra such that the order of lattice is total.\medskip \qquad

\textbf{Definition 5. }([CHA; 58])\textbf{\ }An \textit{MV-algebra }is an
algebra $\left( L,\oplus ,^{\ast },0\right) $ satisfying the following
axioms:

(1) $\left( L,\oplus ,0\right) $ \ is an abelian monoid;

(2) $(x^{\ast })^{\ast }=x;$

(3) $x\oplus 0^{\ast }=0^{\ast };$

(4) $\left( x^{\ast }\oplus y\right) ^{\ast }\oplus y=$ $\left( y^{\ast
}\oplus x\right) ^{\ast }\oplus x$, for all $x,y\in L.\medskip $

\textbf{Remark 6.} If in a BL- algebra $L$ we have $x^{\ast \ast }=x,$ for
every $x\in L$, and, we denote 
\begin{equation*}
x\oplus y=(x^{\ast }\odot y^{\ast })^{\ast },\text{for }x,y\in L,
\end{equation*}%
then we obtain an MV-algebra structure $(L,\oplus ,^{\ast },0).$ Conversely,
if $(L,\oplus ,^{\ast },0)$ is an MV\textit{-}algebra, then $(L,\wedge ,\vee
,\odot ,\rightarrow ,0,1)$ is a BL-algebra, with the following operations: 
\begin{equation*}
x\odot y=(x^{\ast }\oplus y^{\ast })^{\ast },
\end{equation*}%
\begin{equation*}
x\rightarrow y=x^{\ast }\oplus y,1=0^{\ast },
\end{equation*}%
\begin{equation*}
x\vee y=(x\rightarrow y)\rightarrow y=(y\rightarrow x)\rightarrow x\text{
and }x\wedge y=(x^{\ast }\vee y^{\ast })^{\ast },\text{for }x,y\in L.
\end{equation*}%
(see [T; 99]).\medskip

\begin{equation*}
\end{equation*}

\textbf{2. Connections between some polynomial rings and MV-algebras}%
\begin{equation*}
\end{equation*}

From the above Definition 5, we remark that an MV-algebra $\left( L,\oplus
,^{\ast },0\right) $ satisfies some axioms, one of them, $(x^{\ast })^{\ast
}=x$, for all $x\in L,$ attracted our attention in the sense that this
property can be used in defining some new cryptosystems. Ideea behind this
new approach was given by the NTRU cryptosystem, which is a public key
cryptosystem(PKC), where the polynomials are used in defining the public and
the secrete keys. Details about of NTRU cryptosystem and some of its
applications can be found in [TT; 17]. In [CFDP; 22], was proved that if $R$
is a ring factor of a principal integral domain, therefore $(Id(R),\cap
,+,Ann,0=\{0\},1=R)$ is an MV-algebra. To present our cryptosystem, wich is
not PKC, we will use special types of finite principal ideal rings and all
MV-algebras are finite.\medskip\ 

\textbf{Proposition 7. }([CFDP; 22]) \textit{If} $K$ \textit{is a field and} 
$f\in K[x]$ \textit{a polynomial,} $R=K\left[ x\right] /\left( f\right) $%
\textit{, the quotient ring, then} $Id\left( R\right) $ \textit{is an
MV-algebra.}$\Box \medskip $

In the following, we will consider the principal ideal ring $\mathcal{R}%
_{p,1,\beta }=K[x]/\left( x\left( 1-x^{\beta }\right) \right) $. Let $K=%
\mathbb{Z}_{p}$ and $\chi _{\beta }\left( x\right) =x^{\beta +1}-x$. The
lattice $Id\left( \mathcal{R}_{p,1,\beta }\right) $ is an MV-algebra with $%
I^{\ast }=Ann\left( I\right) $ and $I^{\ast \ast }=I,$ for all $I\in
Id\left( \mathcal{R}_{p,1,\beta }\right) $.\medskip

\textbf{Proposition 8.} \textit{Let} $f\in \mathbb{Z}_{p}\left[ x\right] $, $%
1\leq ~$\textit{deg}$\left( f\right) \leq \beta ,$ \textit{such that} $%
f^{2}=1$ \textit{in} $\mathcal{R}_{p,1,\beta }=\mathbb{Z}_{p}[x]/\left(
x\left( 1-x^{\beta }\right) \right) $, \textit{that means} $f=f^{-1}$. 
\textit{Then, there is a natural number} $\delta $ \textit{such that} $f\neq
f^{-1}$ \textit{in} $\mathcal{R}_{p,1,\delta }=\mathbb{Z}_{p}[x]/\left(
x\left( 1-x^{\delta }\right) \right) $.\medskip

\textbf{Proof.} Supposing that that $f^{2}=1$ in $\mathcal{R}_{p,1,\beta }=%
\mathbb{Z}_{p}[x]/\left( x\left( 1-x^{\beta }\right) \right) $, then there
is a polynomial $g\in \mathbb{Z}_{p}[x]$ such that $f\left( x\right) ^{2}+$ $%
g\left( x\right) \left( x^{\beta +1}-x\right) =1$, by using the Euclidean
algorithm. From here, we obtain that $f\left( x\right) ^{2}+$ $g\left(
x\right) x\left( x^{\beta }-1\right) =1$, therefore $f\left( x\right)
^{2}\left( x^{\beta }+1\right) +$ $g\left( x\right) x\left( x^{\beta
}-1\right) \left( x^{\beta }+1\right) =x^{\beta }+1$. It results 
\begin{equation}
f\left( x\right) \rho \left( x\right) +g\left( x\right) \chi _{2\beta
}\left( x\right) =x^{\beta }+1,  \tag{1}
\end{equation}%
where $\rho \left( x\right) =f\left( x\right) \left( x^{\beta }+1\right) $.
Since \textit{deg}$\left( g\right) <\beta $, it is clear that $x^{\beta }+1$
can't be a divisor for $g\left( x\right) $, then relation $\left( 1\right) $
can't have the form $f\left( x\right) ^{2}+$ $g^{\prime }\left( x\right)
\chi _{2\beta }\left( x\right) =1$, where $g\left( x\right) =\left( x^{\beta
}+1\right) g^{\prime }\left( x\right) $. From here, we deduce that the
inverse of the polynomial $f$, if it exists, is different from $f$ in $%
\mathcal{R}_{p,1,2\beta }$, therefore $\delta =2\beta $.$\Box \medskip $

\textbf{Remark 9.} It is obviously that the polynomial $\chi _{p-1}\left(
x\right) =x^{p}-x\in \mathbb{Z}_{p}\left[ x\right] $ has the following
factor decomposition over $\mathbb{Z}_{p}$: $\chi _{p-1}\left( x\right)
=x\left( x+1\right) \left( x-1\right) \left( x+2\right) \left( x-2\right)
...\left( x-\frac{p-1}{2}\right) \left( x+\frac{p-1}{2}\right) $.\medskip

\textbf{The Algorithm. }Let $\mathbb{A}$ be an alphabet with $\lambda $
letters and $M$ a message of length $l$ to be encrypted. The message $M$
received a number $m$ formed by the labels of the componend letters, one by
one, not in blocks. This number is wrote in decimals.

-We consider $p$ a prime number and the polynomial $\chi _{p-1}\left(
x\right) =x^{p}-x$. We convert $m$ in base $p$ and we obtain $m_{p}=%
\overline{a_{q}a_{q-1}...a_{1}}$, with $q\leq p$, $a_{1},a_{2},...,a_{q}\in 
\mathbb{Z}_{p}$. We consider \textit{the associated polynomial message} $%
f_{c}=a_{q}x^{q-1}+a_{q-1}x^{q-2}+...+a_{1}\in \mathbb{Z}_{p}[x]$.

-We consider the field $\mathcal{R}_{p,1,\beta }=\mathbb{Z}_{p}[x]/\left(
x\left( 1-x^{\beta }\right) \right) $, wich is a principal ideal ring, and
we compute its proper ideals, $I_{1},I_{2,}...,I_{j}$. Let $I_{s}=\left(
g_{s}\right) $, where $g_{s}$ is the generator of the Ideal $I_{s}$.

-We found the ideal $I_{t},t\leq j,$ such that $f_{c}\in I_{t},$that means $%
f_{c}\left( x\right) =g_{t}\left( x\right) h\left( x\right) $.

-We compute $Ann\left( I_{t}\right) =I_{r}=\left( g_{r}\right) $ and we
consider \textit{the encrypted polynomial message} $\overline{f_{e}}\left(
x\right) =g_{r}\left( x\right) h\left( x\right)
=b_{v}x^{v-1}+b_{v-1}x^{v-2}+...+b_{1}\in \mathbb{Z}_{p}[x]$. Let $c_{p}=%
\overline{b_{v}b_{v-1}...a_{1}}$ the number in base $p$, which is $c$ in
decimals. We convert $c$ in letters and we get $C$ the encrypted message.

-Since the ideals of the ring $\mathcal{R}_{p,1,\beta }$ form an MV-algebra,
we have that $Ann\left( Ann(I)\right) =I$, that means $Ann\left(
I_{t}\right) =I_{r}$ and $Ann\left( I_{r}\right) =I_{t}$. This remark allows
us decryption of the message, as the rverse of the above steps. The secret
key is $\mathcal{K}=\left( p,\beta ,l\right) $, $p$ a prime numbers, $\beta
+1$ the degree of the polynomial $\chi _{\beta }\left( x\right) $, $\beta $
or $\beta +1$ not necessary to be prime numbers, $l~$the length of the
message. For the situation when the decrypted message has length $l-1$, that
means the message starts with \textbf{A }and\textbf{\ }this implies
insertion of $0$ on the first position in $m$.\medskip

\textbf{Remark 10.} 1) In the ring $\mathcal{R}_{p,1,\beta }$ elements are
invertible or zero divisors. If we obtain that the \textit{\ }attached
polynomial message $f_{c}$ is invertible in $\mathcal{R}_{p,1,\beta }$ and
its inverse, $f_{c}^{-1},$ is different from $f_{c}$, then $f_{c}^{-1}$,
obtained with the extended Euclid's algorithm, is the encrypted polynomial
message $\overline{f_{e}}$. If $f_{c}=f_{c}^{-1}$, then applying Proposition
7, we can find a number $\delta $ such that $f\neq f^{-1}$ in $\mathcal{R}%
_{p,1,\delta }=\mathbb{Z}_{p}[x]/\left( x\left( 1-x^{\delta }\right) \right) 
$ and we apply the algorithm in the ring $\mathcal{R}_{p,1,\delta }$.

2) Usually, $\beta +1\neq p$, but if we take $\beta +1=p$, we can use the
Remark 8, and the ideals of the ring $\mathcal{R}_{p,1,\beta }$ can easily
be computed.\medskip

\textbf{Complexity of the Algorithm.} 1) For the ring $\mathcal{R}%
_{p,1,\beta }=\mathbb{Z}_{p}[x]/\left( x\left( 1-x^{\beta }\right) \right) $%
. In this case, the complexity of this algorithm is influenced by the
multiplication of two polynomials, factors decomposition of a plynomial,
converting a number from decimals to a base $a$ and vice-versa, and the
extended Euclid's algorithms. Multiplication and division of two polinomials
has $O\left( n\log n\right) $ complexity, with $n$ the maximum degree of
those polynomials; \ extended Euclid's algorithm has $O\left( n\left( \log
n\right) ^{2}\right) $; to find an inverse the complexity is $O\left(
n^{2}\log n\log p\right) ,p$ the characteristic of the finite field; to
convert a number $N$ to a base $a,$the complexity is $O\left( N\right) $.
Since the factorization of the polynomial $\chi _{\beta -1}\left( x\right)
=x^{\beta }-x$ is easy to be obtained over $\mathbb{Z}_{p},$ therefore, the
complexity of this algorithm is $O\left( Nn^{2}\left( \log n\right) ^{2}\log
p\right) $.

2) We intend to extend this algorithm, in a further research, to a
commutative principal Artinian ring, as for example is the ring $R=K\left[ x%
\right] /\left( f\right) ,K$ a finite field, $f$ a polynomial of degree $m$,
as we can see in the below next remark. In this case, the above complexity
is influenced by the factoring a polynomial $f$ of degre $m$, such an
algorithm having complexity $O\left( m^{3/2}\log p+m\log ^{2}p\right) $.
Therefore, with the above notations, in this case, the complexity of such an
encryption algorithm is $O\left( Nn^{3}\log ^{2}n\log ^{2}p\left( 1+\log
p\right) \right) $.\medskip

\textbf{Remark 11.} Let $R$ be a commutative, principal, Artinian ring and $%
I\subset R$ an ideal. Therefore $Ann\left( Ann\left( I\right) \right) =I$.
Indeed, since an Artinian ring is finite direct product of Artinian local
rings, then we consider $R$ local. Let $M$ be the unique maximal ideal in $R$%
. If $x\in R$, then $x\in M$ or $x$ is a unit, since in this situation the
set of nonunits form the maximal ideal ~$M$. Ideal $M$ is nilpotent, due the
propery of descending chain of ideals in an Artinian ring, therefore, there
is $t$ such that $M^{t}=\left( 0\right) $. Let $x\in M$ a nonzero element
and $M=\left( x\right) ,$ since the ring is principal. Let $I$ be a nonzero
ideal and $a\in M$ such that $\left( a\right) =I\subseteq M$. We prove that
there is a $k$ such thet $\left( a\right) =M^{k}$. It is clear that $k$ is
such that $a\in M^{k}-M^{k+1}$, since $\left( 0\right) =M^{t}\subseteq
...\subseteq M^{k}\subseteq M^{k-1}\subseteq ...\subseteq M$ $\subseteq R$
is a decreasing sequence. Since $a\in M^{k}$, then $\left( a\right)
\subseteq M^{k}$ and $\widehat{a}\in M^{k}/M^{k-1}$ is nonzero and $%
M^{k}/M^{k-1}$ has dimension $1$, as a vector space, over the field $R/M,$
therefore $M^{k}=\left( a\right) $ and $a=ux^{k}$, $u$ a unit. Therefore $%
I=M^{k}$ and $Ann\left( I\right) =M^{t-k}$. It results, $Ann\left( Ann\left(
I\right) \right) =M^{k}=I$. We obtain that the lattice of ideals of a
commutative, principal, Artinian ring is an MV-algebra. As a general case,
we can take all rings which are are direct sums of local Artinian chain
rings with unit.\medskip\ \medskip

\textbf{Example 12.} 1) If we take $K=\mathbb{Z}_{3},p=3$ and $\beta =2,$
therefore the polynomial $\chi _{2}\left( x\right) =x^{3}-x=$ has the
following decomposition: $x\left( x+1\right) \left( x-1\right) =x\left(
x+1\right) \left( x+2\right) \in \mathbb{Z}_{3}\left[ x\right] $. To avoid a
longue calculus, we consider an alphabet with $10$ letters, labeled as in
the below table: 
\begin{equation*}
\begin{tabular}{|l|l|l|l|l|l|l|l|l|l|}
\hline
$A$ & $B$ & $C$ & $D$ & $E$ & $F$ & $G$ & $H$ & $I$ & $J$ \\ \hline
$0$ & $1$ & $2$ & $3$ & $4$ & $5$ & $6$ & $7$ & $8$ & $9$ \\ \hline
\end{tabular}%
.
\end{equation*}%
The ideals of the ring $\mathcal{R}_{3,1,2}~$are: $\left( 0\right) ,\mathcal{%
R}_{3,1,2},\left( x\right) ,\left( x-1\right) ,\left( x+1\right) ,\left(
x^{2}-x\right) ,\left( x^{2}+x\right) ,\left( x^{2}-1\right) $, in total, $8$
ideals. We want to encrypt the message \textbf{BJ}. Its decimal label is $%
m=19$, which is $m_{3}=201$ in base $3$. The associated polynomial is $%
f_{c}\left( x\right) =2x^{2}+1=2\left( x+1\right) \left( x-1\right) =2\left(
x^{2}-1\right) \in I_{t}=\left( x^{2}-1\right) $. We have $f_{c}\left(
x\right) =g_{t}\left( x\right) h\left( x\right) =2\left( x^{2}-1\right)
,h\left( x\right) =2$ and $Ann\left( I_{t}\right) =\left( x\right) $,
therefore the encrypted polynomial\textit{\ }message $\overline{f_{e}}\left(
x\right) =2x$. We obtain $c_{3}=$ $020$ in base $3$ wich is $c=6$ in
decimal. Therefore, the encrypted message is \textbf{G}. In this case, the
encryption key is $\mathcal{K}=\left( 3,2,2\right) $.

2) We take $K=\mathbb{Z}_{3},p=3$ and $\beta =4,$ therefore the polynomial $%
\chi _{4}\left( x\right) =x^{5}-x$ has the following decomposition: $%
x^{5}-x=x\left( x-1\right) \left( x+1\right) \left( x^{2}+1\right) \in 
\mathbb{Z}_{3}\left[ x\right] $. We want to encrypt the message \textbf{ABBA}%
. The attached decimal label is $m=0110$, which is $m_{3}=11002$ in base $3$%
. The key in this situation is $\mathcal{K}=\left( 3,4,4\right) $. The
associated polynomial is $f_{c}=x^{4}+x^{3}+2$, which is an invertible
element in $\mathcal{R}_{3,1,4},$ with $f_{c}^{\prime }=$ $x^{4}+x+2$ its
inverse. The label will be $c_{3}=10012$, in base $3$, which is $c=84$ in
decimal, therefore the encrypted text is \textbf{IE}. If we want decrypt
this message, we find $\left( 84\right) _{3}=10012$, the attached polynomial
is $f_{c}^{\prime }$, with its inverse $f_{c}$, and we obtain $c=110$ in
decimals. Since from the transmitted key, the length of the message is $4\,$%
, this imply that we have a $0$ on the first position, then $0110\rightarrow
ABBA$, is the decrypted message.

3) The above message \textbf{ABBA}, can be encrypted in another way, namely
if we consider $p=5,$then the encryption key is $\mathcal{K}=\left(
5,2,4\right) $. Therefore, we have $K=\mathbb{Z}_{5},p=5$ and $\beta =2,$
and the polynomial $\chi _{2}\left( x\right) =x^{3}-x$ has the following
decomposition: $x\left( x+1\right) \left( x+4\right) \in \mathbb{Z}_{5}\left[
x\right] $. The attached decimal label $m=0110$, which is $m_{5}=420$ in
base $5$ and the associated polynomial is $f_{c}=4x^{2}+2x=x\left(
4x+2\right) \in \left( x\right) $. Since $Ann\left( \left( x\right) \right)
=(\left( x+1\right) \left( x+4\right) )=\left( \left( x^{2}-1\right) \right)
,$ the encrypted polynomial\textit{\ }message will be $\overline{f_{e}}%
\left( x\right) =$ $\left( x^{2}-1\right) \left( 4x+2\right) =2x^{2}+3$.
Then, the label is $c_{5}=203$ in base $5$, which is $c=53$ in decimal. The
encrypted text is \textbf{FD}. To decrypt the message \textbf{FD}, $53$
becomes $203$ in base $5$, with the associated polynomial $2x^{2}+3\in
\left( x^{2}-1\right) $, with the quotient polynomial $q\left( x\right)
=\left( 4x+2\right) $. We have $Ann\left( \left( x^{2}-1\right) \right)
=\left( x\right) $, then the decryption polynomial is $d\left( x\right)
=\gamma \left( x\right) x=4x^{2}+2x$, which give us the label $420$ in base $%
5$. We obtain $110$ in decimal, then \textbf{BBA}. Since the length of the
message is $4$, we have a $0$ on the first position, then $0110\rightarrow
ABBA$ is the decrypted message.

4) We take $K=\mathbb{Z}_{3},p=3$ and $\beta =2,$ therefore $\mathcal{R}%
_{3,1,2}=\mathbb{Z}_{3}[x]/\left( x\left( 1-x^{2}\right) \right) $. The
plain text is \textbf{CF}, with decimal label $m=25$ and $m_{3}=221$ in base 
$3$. The associated polynomial is $f_{c}\left( x\right) =2x^{2}+2x+1$, with $%
f^{2}=1$ in $\mathcal{R}_{3,1,2}$ and $f^{-1}=f$. Therefore, we consider the
ring $\mathcal{R}_{3,1,4}=\mathbb{Z}_{3}[x]/\left( x\left( 1-x^{4}\right)
\right) $ and $f^{-1}=x^{4}+x^{2}+2x+1$. The obtained label in base $3$ is $%
c_{3}=10121$. In decimal base will be $c=97$, then the encrypted message is 
\textbf{JH}. The secret key is $\left( 3,4,2\right) $.

5) We want encrypt the text \textbf{DECADE}. We obtain $m=342034$ and $%
m_{3}=122101011221$, in base $3,m_{5}=41421114$, in base $5$ and $%
m_{7}=2623120$, in base $7$. Since $m_{7}$ has the smaller length, we will
consider $p=7,\mathcal{R}_{7,1,6}=\mathbb{Z}_{3}[x]/\left( x\left(
1-x^{6}\right) \right) $ and $\chi _{p-1}\left( x\right) =x^{7}-x=x\left(
x+1\right) \left( x+6\right) \left( x+2\right) \left( x+5\right) \left(
x+3\right) \left( x+4\right) $. In this situation, the encryption key is $%
\mathcal{K}=\left( 7,6,6\right) $. The associated polynomial message $f_{c}$
is $f_{c}=2x^{6}+6x^{5}+2x^{4}+3x^{3}+x^{2}+2x=x\left( x+2\right) \left(
2x^{4}+2x^{3}+5x^{2}+1\right) \in \left( x\left( x+2\right) \right) $, where 
$I_{t}=\left( x\left( x+2\right) \right) $ is the ideal generated by the
polynomial $g_{t}\left( x\right) =x^{2}+2x$ and $h\left( x\right)
=2x^{4}+2x^{3}+5x^{2}+1$. The $Ann\left( I_{t}\right) =I_{r}=\left(
g_{r}\right) $, $g_{r}\left( x\right) =\left( x+1\right) \left( x+6\right)
\left( x+5\right) \left( x+3\right) \left( x+4\right) $.

We obtain the encrypted polynomial message $\overline{f_{e}}\left( x\right)
=g_{r}\left( x\right) h\left( x\right)
=3x^{6}+2x^{5}+3x^{4}+x^{3}+5x^{2}+4x+3$ and $c_{7}=3231543$. In decimals, $%
c_{7}$ is $c=394383$ and the encrypted message is \textbf{DJEDID}.

\begin{equation*}
\end{equation*}

\textbf{3. Remarks regarding BL-comets}%
\begin{equation*}
\end{equation*}

In the paper [NL; 03], authors analyzed the structure of finite BL-algebras.
They introduced the concept of BL-comets, a class of finite BL-algebras
which can be seen as a generalization of finite BL-chains. Using BL-comets,
any finite BL-algebra can be representd as a direct product of
BL-comets.\medskip\ 

\textbf{Definition 13.}( [NL; 05], Definition 3, [FP; 22]) Let $\left(
C_{i},\wedge _{i},\vee _{i},\odot _{i},\rightarrow _{i},0_{i},1_{i}\right)
,i\in \{1,2,...,t-1\}$ be $t-1$ BL-chains and $C_{t}$ a BL-algebra. We
consider $1_{i}=0_{i+1},i\in \{1,2,...,t-1\},0=0_{1}$, $1=1_{t}$ and that $%
(C_{i}\backslash \{1_{i}\})\cap (C_{i+1}\backslash \{0_{i+1}\})=\emptyset ,$%
for $i\in \{1,2,...,t-1\}$. \ The \textit{ordinal sum} $\underset{i=1}{%
\overset{t}{\uplus }}C_{i}$ is defined to be the following BL-algebra 
\begin{equation*}
\left( \overset{t}{\underset{i=1}{\cup }}C_{i},\wedge ,\vee ,\odot
,\rightarrow ,0,1\right) ,
\end{equation*}%
whose operations are defined as follows

\begin{equation*}
x\leq y\text{ if }(x,y\in C_{i}\text{ and }x\leq _{i}y)\text{ or }(x\in C_{i}%
\text{ and }y\in C_{j},i<j,i,j\in \{1,2,...,t\})\text{ ,}
\end{equation*}

\begin{equation*}
x\wedge y=\left\{ 
\begin{array}{c}
x\wedge _{i}y,\text{ if }x,y\in C_{i}, \\ 
x,\text{ if }x\in C_{i}\text{ and }y\in C_{j},i<j,i,j\in \{1,2,...,t\}%
\end{array}%
\right.
\end{equation*}

\begin{equation*}
x\vee y=\left\{ 
\begin{array}{c}
x\vee _{i}y,\text{ if }x,y\in C_{i}, \\ 
y,\text{ if }x\in C_{i}\text{ and }y\in C_{j},i<j,i,j\in \{1,2,...,t\}%
\end{array}%
\right.
\end{equation*}%
\begin{equation*}
x\rightarrow y=\left\{ 
\begin{array}{c}
1,\text{ if }x\leq y, \\ 
x\rightarrow _{i}y,\text{ if }x\nleq y,\text{ }x,y\in C_{i},\text{ }i\in
\{1,2,...,t\}, \\ 
y,\text{ if }x\nleq y,\text{ }x\in C_{j},\text{ }y\in C_{i}\backslash
\{1_{i}\},i<j.%
\end{array}%
\right.
\end{equation*}%
\begin{equation*}
x\odot y=\left\{ 
\begin{array}{c}
x\odot _{i}y,\text{ if }x,y\in C_{i}, \\ 
x,\text{ if }x\in C_{i}\backslash \{1_{i}\}\text{ and }y\in C_{j},i<j%
\end{array}%
\right. .
\end{equation*}%
We will write $\underset{i=1}{\overset{t}{\uplus }}C_{i}$ as $C_{1}\boxplus
C_{2}\boxplus ...\boxplus C_{t}.$

\textbf{Definition 14.} 1) ([NL; 03], Definition 21) Let $L$ be a
BL-algebra. The element $x\in L$ is called \textit{idempotent} if $x\odot
x=x $.

2) We consider $L$ a finite BL-algebra and $\mathcal{I}\left( L\right) $ the
set of idempotent elements in $L$. For $x\in \mathcal{I}\left( L\right) $,
we denote $\mathcal{C}\left( x\right) =\{y\in \mathcal{I}\left( L\right) ~$\
such that $x$ and $y$ are comparable$\}$. We define the set $\mathcal{D}%
\left( L\right) \subseteq \mathcal{I}\left( L\right) $ as follows:

$x\in \mathcal{D}\left( L\right) $ if and only if

i) $\mathcal{C}\left( x\right) =\mathcal{I}\left( L\right) ;$

ii) The set $\{y\in \mathcal{I}\left( L\right) ,$ $y\leq x\}$ is a chain.

We obtain that $\mathcal{D}\left( L\right) \neq \emptyset ,$ since $0\in 
\mathcal{D}\left( L\right) .$

A finite BL-algebra $L$ is called a \textit{BL-comet} if \textit{max}$%
\mathcal{D}\left( L\right) \neq 0.$

In a BL-comet $L,$ the element \textit{max}$\mathcal{D}\left( L\right) $ is
called \textit{the pivot} of $L$ and it is denoted by $pivot(L).\medskip $

\textbf{Proposition 15. }([NL; 03], Proposition 26) \textit{Let} $L$ \textit{%
be a finite BL-algebra. The following assertions are equivalent:}

(i) $L$ \textit{is a BL-comet and }$pivot(L)=1;$

(ii)\textit{\ }$L$ \textit{is a BL-chain}.\medskip $\Box $

\textbf{Remark 16.} 1) From [NL, 03], a finite BL-chain is defined to be a
finite ordinal sum of finite MV-chains. In the same paper, authors analyzed
the structure of finite BL-algebras and introduced the concept of BL-comets,
a class of finite BL-algebras which can be seen as a generalization of
finite BL-chains. Using BL-comets, they proved that any finite BL-algebra
can be represent as a direct product of BL-comets (Corollary 10). From here,
we have that a finite BL-algebra $L$ with a prime number of elements is a BL
chain or a comet with $pivot(L)<1$

2) ([I; 09], Corollary 3.5.10) If $L_{1}$ and $L_{2}$ are two BL-algebras
and $L_{1}$ is a BL-chain, then the ordinal sum $L_{1}\boxplus L_{2}$ is a
BL-algebra.\medskip

\textbf{Proposition 17.} ([NL; 05], Theorem 22 and Corollary 24) \textit{Let}
$L$ \textit{be a finite BL-algebra. If }$L$ \textit{is a BL-comet with }$%
pivot(L)<1$\textit{, then} $L$ \textit{is the ordinal sum of a finite
BL-chain and a finite BL-algebra which is not a BL-comet.\medskip }$\Box $

\textbf{Proposition 18. (}[CFP; 23])

1) \textit{Let} $L$ \textit{be a BL-comet. Then }$L$ \textit{is a BL-chain
iff} $pivot(L)^{\ast \ast }=pivot(L).$

2) \textit{Let} $L$ \textit{be a finite MV-algebra. The following assertions
are equivalent:}

(i) $L$ \textit{is a} \textit{BL-comet;}

(ii)\textit{\ }$L$ \textit{is an MV-chain.}$\Box \medskip $

The ideea of this section arised from the fact that in our researches we try
to find types of rings $R$ such that on $Id\left( R\right) $, if it is a
finite set, to obtain a BL-algebra structures which are not MV-algebras. But
a commonplace example of order three 
\begin{equation*}
\begin{tabular}{l|lll}
$\rightarrow $ & $0$ & $a$ & $1$ \\ \hline
$0$ & $1$ & $1$ & $1$ \\ 
$a$ & $0$ & $1$ & $1$ \\ 
$1$ & $0$ & $a$ & $1$%
\end{tabular}%
\ \ \ 
\begin{tabular}{l|lll}
$\otimes $ & $0$ & $a$ & $1$ \\ \hline
$0$ & $0$ & $0$ & $0$ \\ 
$a$ & $0$ & $a$ & $a$ \\ 
$1$ & $0$ & $a$ & $1$%
\end{tabular}%
\end{equation*}%
gives us a BL-algebra which is not an MV-algebra, such that there is not a
commutative unitary ring $R$ with three ideals, with the algebra $Id(R)$
being a BL-algebra, with $\rightarrow $ and $\otimes $ defined above. This
is an example of BL-chain wich is non an MV-chain. As we can see, a BL-chain
is a particular case of a BL-comet. We asked if this situation is an isolate
case or can be generalised. Indeed, this result can be extended, to all
BL-comet, chain or not, as we can see in Theorem 31.\medskip

\textbf{Proposition 19.} (see [CFP; 23]) \textit{Let} $R$ \textit{be a
commutative and unitary ring with a finite number of ideals. Let} $%
n_{m}\left( R\right) $ \textit{be the number of maximal ideals in} $R$%
\textit{,} $n_{p}\left( R\right) $ \textit{be the number of prime ideals in} 
$R$ \textit{and} $n_{I}\left( R\right) $ \textit{be the number of all ideals
in} $R$\textit{. Therefore,} $n_{m}\left( R\right) =n_{p}\left( R\right)
=\alpha $ \textit{and} $n_{I}\left( R\right) =\underset{j=1}{\overset{\alpha 
}{\prod }}\beta _{j},\beta _{j}$ \textit{positive integers, }$\beta _{j}\geq
2$. $\Box \medskip $

\textbf{Example 20.} In [FP; 22], we presented a basic summary of the
structure of BL-algebras with $n$ elements, $2\leq n\leq 5$. For $n=5$, were
obtained $9$ different types, namely: 
\begin{equation*}
\left\{ 
\begin{array}{c}
Id(\mathbb{Z}_{16})\text{ (chain, MV)} \\ 
Id(\mathbb{Z}_{2})\boxplus Id(\mathbb{Z}_{8})\text{ (BL-chain)} \\ 
Id(\mathbb{Z}_{2})\boxplus Id(\mathbb{Z}_{2}\times \mathbb{Z}_{2})\text{
(comet)} \\ 
Id(\mathbb{Z}_{2})\boxplus (Id(\mathbb{Z}_{2})\boxplus Id(\mathbb{Z}_{4}))%
\text{ (BL-chain)} \\ 
Id(\mathbb{Z}_{2})\boxplus (Id(\mathbb{Z}_{4})\boxplus Id(\mathbb{Z}_{2}))%
\text{ (BL-chain)} \\ 
Id(\mathbb{Z}_{2})\boxplus (Id(\mathbb{Z}_{2})\boxplus (Id(\mathbb{Z}%
_{2})\boxplus Id(\mathbb{Z}_{2})))\text{ (BL-chain)} \\ 
Id(\mathbb{Z}_{8})\boxplus Id(\mathbb{Z}_{2})\text{ (BL-chain)} \\ 
(Id(\mathbb{Z}_{4})\boxplus Id(\mathbb{Z}_{2}))\boxplus Id(\mathbb{Z}_{2})%
\text{ (BL-chain)} \\ 
Id(\mathbb{Z}_{4})\boxplus Id(\mathbb{Z}_{4})\text{ (BL-chain)}%
\end{array}%
\right. .
\end{equation*}%
The lattice $\mathcal{L}_{5}=Id(\mathbb{Z}_{2})\boxplus Id(\mathbb{Z}%
_{2}\times \mathbb{Z}_{2})$ is a BL-comet lattice. Indeed, this lattice $%
\mathcal{L}_{5}=\{0,a,b,c,1\}$ is a finite BL-algebra which is not an
MV-algebra and has the following operations: 
\begin{equation*}
\begin{array}{c|ccccc}
\rightarrow & 0 & a & b & c & 1 \\ \hline
0 & 1 & 1 & 1 & 1 & 1 \\ 
a & 0 & 1 & 1 & 1 & 1 \\ 
b & 0 & c & 1 & c & 1 \\ 
c & 0 & b & b & 1 & 1 \\ 
1 & 0 & a & b & c & 1%
\end{array}%
,\hspace{5mm}%
\begin{array}{c|ccccc}
\odot & 0 & a & b & c & 1 \\ \hline
0 & 0 & 0 & 0 & 0 & 0 \\ 
a & 0 & a & a & a & a \\ 
b & 0 & a & b & a & b \\ 
c & 0 & a & a & c & c \\ 
1 & 0 & a & b & c & 1%
\end{array}%
,
\end{equation*}%
where $Id(\mathbb{Z}_{2})=\{0,a\},a=\mathbb{Z}_{2},0=\left( 0\right) $ and $%
Id(\mathbb{Z}_{2}\times \mathbb{Z}_{2})=\{0,b,c,1\},$with $0=\left( 0\right)
,1=\mathbb{Z}_{2}\times \mathbb{Z}_{2},b=\{\left( 0,0\right) ,\left(
0,1\right) \},c=\{\left( 0,0\right) ,\left( 1,0\right) \}$. We have that $%
\mathcal{L}_{5}$ is a BL-comet. Indeed, by using definition of a BL-comet,
we have $\mathcal{I}\left( L_{5}\right) =\{0,a,b,c,1\}$. We take $x=a$, then 
$\mathcal{C}\left( a\right) =$ $\mathcal{I}\left( \mathcal{L}_{5}\right) $
and the set $\{y\in \mathcal{I}\left( \mathcal{L}_{5}\right) ,y\leq
a\}=\{0,a\}$ is a chain. Therefore, $\mathcal{D}\left( \mathcal{L}%
_{5}\right) =\{0,a\}$ with $pivot=$~$\ $\textit{max}$\mathcal{D}\left( 
\mathcal{L}_{5}\right) =a\neq 0$, $a<1$. Since $a<1$, we have that $\mathcal{%
L}_{5}$ is the ordinal sum of a finite BL-chain and a finite BL-algebra
which is not a BL-comet: $Id(\mathbb{Z}_{2})$ is a BL-chain and $Id(\mathbb{Z%
}_{2}\times \mathbb{Z}_{2})$ is an MV-algebra (BL) wich is not a BL-comet.
We remark that $\mathcal{L}_{5}$ has two maximal elements, $b$ and $c$,
which correspond to the two maximal ideals of the ring $\mathbb{Z}_{2}\times 
\mathbb{Z}_{2}$.\medskip

\textbf{Definition 21.} Let $L$ be a BL-algebra and $x,y\in L$. We have that 
$x\leq y$ iff $x\rightarrow y=1$. The element $m\in L$ is called a \textit{%
maximal element} in $L$ if and only if for each $x\in L$ such that $x\leq m$%
, we have $x\rightarrow m=1$ and if $m\leq y$, we have $m=y$ or $y=1$. The
dual concept of a maximal element in $L$ is the \textit{minimal element}%
.\medskip

\textbf{Remark 22.} If $L$ is a BL-algebra such that there is a ring $R$
with $Id\left( R\right) =L$, then maximal ideals in $R$ are maximal elements
in $L$ and vice-versa and the minimal ideals in $R$ are minimal elements in $%
L$ and vice-versa.\medskip

\textbf{Proposition 23.} ([CFP]) \textit{Let} $R$ \textit{be a commutative
unitary ring which has exactly three ideals }$\{0\},I,R$. \textit{Therefore,
we have} $I^{2}=\{0\}$.

ii) \textit{There are no commutative unitary rings} $R$ \textit{with three
ideals having} $(Id(R),\cap ,+,\otimes \rightarrow ,0=\{0\},1=R)$ \textit{as
a BL-algebra which is not an MV-algebra.}$\Box $\textit{\medskip }

\textbf{Proposition 24.} \textit{A local ring} $R$ \textit{doesn't contains
nontrivial idempotents}.\medskip\ 

\textbf{Proof}. Indeed, if \ $e$ is an idempotent, $e\neq 0,1$, then $%
e\left( e-1\right) =e^{2}-e=0$. From here, we have that $e$ and $e-1$ are
non-invertible zero-divisors and belong to the unique maximal ideal $M$.
Since $1=e+\left( 1-e\right) $, we obtain that $1\in M$, then $M=R$, false. $%
\Box \medskip $

\textbf{Proposition 25. }\textit{If} $L$ \textit{is a BL-comet, with }$%
pivot\left( L\right) =1$\textit{(that means a BL-chain), then there are no
commutative and unitary rings} $R$ \textit{such that} $Id\left( R\right) =L$%
.\medskip\ 

\textbf{Proof.} From the above, we have that $L$ is a BL-chain and it is a
finite ordinal sum of finite MV-chain, $\left( M_{i},0_{i},1_{i}\right)
,i\in \{1,2,...,t\}$, $L=\underset{i=1}{\overset{t}{\uplus }}M_{i}$. For $%
i\in \{1,2,...,t-1\}$, the element $a_{i}=1_{i}=0_{i+1}$ is an nontrivial
idempotent in $L$. If there is a ring $R$ such that $Id\left( R\right) =L$,
then $R$ is a local ring and hasn't nontrivial idempotents, false.$\Box
\medskip $

\textbf{Proposition 26.} \textit{We consider} $L$ \textit{a finite BL-comet
algebra, with} $\left\vert L\right\vert =n$. \textit{If} $L$ \textit{is a
BL-chain, then} $L$ \textit{has only one maximal element and only one
minimal element. If} $L$ \textit{is not a chain, then} $L$ \textit{has
minimum two maximal elements and only one minimal element.\medskip }

\textbf{Proof. }If $L$ is a chain, it is clear that has only one maximal
element and only one minimal element. We make induction after $n$.

For $\left\vert L\right\vert =n=2$ and $3$, we have a BL-chain comet,
therefore we have one maximal element and only one minimal element. For $%
\left\vert L\right\vert =4$, $L$ is a BL-chain with one maximal element and
only one minimal element. For $\left\vert L\right\vert =5$, we have that $L$
is a BL-chain with one maximal element and only one minimal element or $L=%
\mathcal{L}_{5}$, as in the above example, and has two maximal elements and
only one minimal element. Assuming that all BL-comets \thinspace $L$, wich
are not chains and $\left\vert L\right\vert <n$, has minimum two maximal
elements and only one minimal element, let $L_{n}$ be a BL-comet with $%
\left\vert L_{n}\right\vert =n$. We have that $L_{n}$ is an ordinal sum
between finite chains $C_{s}~$(then $L_{n}$ has and only one minimal
element) and a finite BL-algebra $B$, $B$ is not comet. Therefore, $B$ is a
direct product of minimum two BL-comets (chain or not), $B=B_{1}\times
...\times B_{t}$, $t\geq 2$, with $\left\vert B_{i}\right\vert <n$. By using
the induction hypothesis, each $B_{i}$ has minimum one maximal element and $%
B $ will have minimum two maximal elements. We remark that, these maximal
elements in $B$ are maximal elements in the BL-comet $L_{n}$, due to the
definition of ordinal sum. We remark that $\left\vert B\right\vert $ is not
a prime number, since in this case $B$ must be a BL-comet, false.$\Box
\medskip $

\textbf{Proposition 27. }\textit{If} $L$ \textit{is a BL-comet, with }$%
pivot\left( L\right) <1$ \textit{and} $\left\vert L\right\vert =p$, $p$ 
\textit{a prime number, then there is no commutative and unitary ring} $R$ 
\textit{such that} $Id\left( R\right) =L$.\medskip\ 

\textbf{Proof.} Supposing that there is a ring $R$ such that $Id\left(
R\right) =L$. From the above proposition, $L$ has at least two maximal
elements, which correspond to two maximal ideals in $R$. Since $\left\vert
Id\left( R\right) \right\vert =n_{I}\left( R\right) =p$, $p$ a prime number,
and $n_{I}\left( R\right) =\underset{j=1}{\overset{\alpha }{\prod }}\beta
_{j},\beta _{j}$ positive integers,\textit{\ }$\beta _{j}\geq 2$, with $%
\alpha =n_{m}\left( R\right) $, the number of maximal ideals, which is at
least two, we have a contradiction.\medskip $\Box $

\textbf{Remark 28.} From the above, we remark that for $n=2,$we have a
chain, for $n=3$,we have an MV-chain, $Id\left( \mathbb{Z}_{4}\right) $ and
a BL-chain, which is not an MV-chain, $Id(\mathbb{Z}_{2})\boxplus Id(\mathbb{%
Z}_{2})=\{\{0\},\{0,1\}\}\boxplus \{\{0\},\{0,1\}\}$, with $%
a=\{0,1\}\boxplus \{0\}$, a nontrivial idempotent element, with the below
multiplication tables: 
\begin{equation*}
\begin{tabular}{l|lll}
$\rightarrow $ & $0$ & $a$ & $1$ \\ \hline
$0$ & $1$ & $1$ & $1$ \\ 
$a$ & $0$ & $1$ & $1$ \\ 
$1$ & $0$ & $a$ & $1$%
\end{tabular}%
\ \ \ 
\begin{tabular}{l|lll}
$\otimes $ & $0$ & $a$ & $1$ \\ \hline
$0$ & $0$ & $0$ & $0$ \\ 
$a$ & $0$ & $a$ & $a$ \\ 
$1$ & $0$ & $a$ & $1$%
\end{tabular}%
.
\end{equation*}

\bigskip Therefore, from the above results, we obtain the following theorem:

\textbf{Remark 29.} 1) (\textbf{[AM; 69],} Proposition 8.1.) In acommutative
unitary Artinian ring $A$, every prime ideal is maximal and vice-versa.

2) We consider $R$ a commutative unitary ring with a finite number of
ideals, which is not a field. The ring $R$ is an Artinian and a Noetherian
ring in the same time. We prove that a prime ideal in the ring $R$ has a
nonzero annihilator, therefore a maximal ideal in such a ring has a nonzero
annihilator. Indeed, let $x\in R\,$\ and $Ann\left( x\right) =\left\{ r\in
R,rx=0\right\} $ be the annihilator of the element $x$. $Ann\left( x\right) $
is an ideal in $R.$ We consider the set 
\begin{equation*}
\mathcal{A}=\left\{ Ann\left( x\right) ,x\in R,x\neq 0\right\} .
\end{equation*}%
It is clear that $\mathcal{A}$ is a finite set, since we have a finite
number of ideals in $R$. Therefore, there is a maximal element in $\mathcal{A%
}$, namely, $J=Ann\left( x\right) $, with $x\neq 0$. The ideal $J$ is a
prime ideal, therefore is a maximal ideal. Indeed, let $\alpha ,\beta \in
R-J $ such that $\alpha \beta \in J$. We have that $\alpha x\neq 0,\beta
x\neq 0$, but $\alpha \beta x=0$, therefore $\alpha \beta \in J=Ann\left(
x\right) $. We consider the set $Ann\left( \alpha x\right) =\left\{ r\in
R,r\left( \alpha x\right) =0\right\} $. It results that $Ann\left( x\right)
\varsubsetneq Ann\left( \alpha x\right) $, with $\alpha x\neq 0$, then $%
Ann\left( \alpha x\right) \in \mathcal{A}$, contradiction with the fact that 
$J$ is the maximal element in $\mathcal{A}$. Therefore, if $\alpha \beta \in
J$, then $\alpha \in J$ or $\beta \in J$ and $J$ is a prime ideal. It
results that $J=Ann\left( x\right) $ is a prime ideal which is the
annihilator of a nonzero element. Therefore, each maximal ideal has a
nonzero annihilator. We remark that if $J=\left( 0\right) $ is prime, this
is equivalent with the fact that $R$ is an integral domain (\textbf{[AM; 69],%
} p. 3) and an integral domain with a finite number of ideals is a filed([%
\textbf{CFP; 23}], Proposition 2.10), contradiction.\medskip

\textbf{Remark 30. }Let $R$ be a commutative and unitary ring with a finite
number of ideals and $M$ a maximal ideal. Since we proved that $Ann\left(
M\right) \neq \left( 0\right) $, then there is a minimal ideal $I_{m}$ such
that $I_{m}\subseteq Ann\left( M\right) $. From here, we have that $I_{m}M=0$%
, then $M\subseteq Ann\left( I_{m}\right) $. Since $M$ is maximal, we have $%
M=Ann\left( I_{m}\right) $. Therefore, for a maximal ideal $M$, always exist
a minimal ideal $I_{m}$ such that $M=Ann\left( I_{m}\right) $.\medskip

\textbf{Theorem 31.} \textit{If} $L$ \textit{is a finite BL-comet, with }$%
pivot\left( L\right) <1$, \textit{then there is no commutative and unitary
rings} $R$ \textit{such that} $Id\left( R\right) =L$.\medskip

\textbf{Proof.} Using results obtained in the above remarks, if there is a
ring $R$ such that $Id\left( R\right) =L$, since $L$ has only one minimal
ideal $J$ and minimum two maximal ideals, $M_{1},M_{2}$, we have that $M_{1}$
and $M_{2}$ are the annulators of some minimal ideals $J_{1},J_{2}$: $%
M_{1}=Ann\left( J_{1}\right) \neq 0$ and $M_{2}=Ann\left( J_{2}\right) \neq
0 $. In our case $J_{1}=J_{2}=J,$therefore $M_{1}=M_{2}$, contradiction.$%
\Box \medskip $

\begin{equation*}
\end{equation*}

\textbf{4. Characterisation of finite BL-algebras}%
\begin{equation*}
\end{equation*}

\textbf{Remark 32.} 1)The ordinal sum of two BL-algebras $\mathcal{L}%
_{1}=(L_{1},\wedge _{1},\vee _{1},\odot _{1},\rightarrow _{1},0_{1},1_{1})$
and $\mathcal{L}_{2}=(L_{2},\wedge _{2},\vee _{2},\odot _{2},\rightarrow
_{2},0_{2},1_{2})$ with $1_{1}=0_{2}$ and $(L_{1}\backslash \{1_{1}\})\cap
(L_{2}\backslash \{0_{2}\})=\emptyset $ \ is a residuated lattice $\mathcal{L%
}_{1}\boxplus \mathcal{L}_{2}=(L_{1}\cup L_{2},\wedge ,\vee ,\odot
,\rightarrow ,0=0_{1}$ $,1=1_{2}$ $)$ which is not a BL algebra if $L_{1}$\
is not a chain. Indeed, if $L_{1}$\ is not a chain, then there are $a,b\in
L_{1}$ incomparable. Then $(a\rightarrow b)\vee (b\rightarrow a)=$ $%
(a\rightarrow _{1}b)\vee (b\rightarrow _{1}a)=1_{1}\neq 1_{2}=1.$

2) The ordinal sum between a BL chain $L_{1}$ and a BL-algebra $L_{2}$ is a
BL-algebra $L_{1}\boxplus L_{2}$ with $\max \mathcal{D}(L_{1}\boxplus
L_{2})\neq 0$ which is not an MV-algebra. Indeed, $L_{1}\boxplus L_{2}$ is a
BL-algebra with 
\begin{equation*}
(1_{1})^{\ast \ast }=(1_{1}\rightarrow 0_{1})\rightarrow
0_{1}=0_{1}\rightarrow 0_{1}=1_{2}\neq 1_{1}.
\end{equation*}

Since $1_{1}=0_{2}\in \mathcal{I}(L_{1}\boxplus L_{2}),$ $\mathcal{C}(1_{1})=%
\mathcal{I}(L_{1}\boxplus L_{2})$ and $\{y\in \mathcal{I}(L_{1}\boxplus
L_{2}):y\leq 1_{1}\}=\{y\in \mathcal{I}(L_{1}):y\leq 1_{1}\}$ is \ a chain,
we deduce that $1_{1}=0_{2}\in \mathcal{D}(L_{1}\boxplus L_{2}),$ so, 
\textit{max} $\mathcal{D}(L_{1}\boxplus L_{2})\neq 0=0_{1.}$

3) Definition 13 provides a way to generate finite\ BL-comets which are not
MV-algebras.\medskip

\textbf{Lemma 33. }\textit{Let} $L$ \textit{be a finite BL-algebra and} $a=~$%
\textit{max}$\mathcal{D}(L).$ \textit{Then} $a=0$ \textit{or} $a^{\ast
}=0.\medskip $

\textbf{Proof.} Obviously, $0\in \mathcal{D}(L).$

Suppose that $a\neq 0.$

We recall that in a BL-algebra $L$, $(x\odot y)^{\ast \ast }=x^{\ast \ast
}\odot y^{\ast \ast },$ for any $x,y\in L.$ For $x=y=a$ we deduce that $%
(a^{2})^{\ast \ast }=(a^{\ast \ast })^{2}.$ Since $a\in \mathcal{I}(L)$ we
deduce that $a^{\ast \ast }=(a^{\ast \ast })^{2},$ so $a^{\ast \ast }\in 
\mathcal{I}(L).$ Using the caracterization of boolean elements in a
BL-algebra (see [P; 07]) we deduce that $a^{\ast \ast }\in \mathcal{B}(L)=$
the set of boolean elements of $L,$ so $a^{\ast }=(a^{\ast \ast })^{\ast
}\in \mathcal{B}(L).$ Then $a^{\ast }\in \mathcal{I}(L).$

Since $\mathcal{C}(a)=\mathcal{I}(L),$ $a$ and $a^{\ast }$ are comparable.

If $a\leq a^{\ast }$ then $0=a\odot a^{\ast }=a\wedge a^{\ast }=a,$ a
contradiction.

If $a^{\ast }\leq a$ then $0=a\odot a^{\ast }=a\wedge a^{\ast }=a^{\ast
}.\Box \medskip $

\textbf{Theorem 34. }\textit{Let} $L$ \textit{be a finite MV-algebra. Then} $%
\max \mathcal{D}(L)\in \{0,1\}.\medskip $

\textbf{Proof 1.} Obviously, from Remark 5, MV-algebras are particular
BL-algebras. Using Proposition 18, an MV-algebra is a chain iff it is a
BL-comet, and for an MV-chain, $\max \mathcal{D}(L)=1.$

If $L$ is not a chain, then obviously, it is not a comet, so $\max \mathcal{D%
}(L)=0.\medskip $

\textbf{Proof 2.} $L$ is in particular a BL-algebra. From Lemma 33, \textbf{%
\ }if $a=\max \mathcal{D}(L),$ then $a=0$ or $a^{\ast }=0.$ If $a\neq 0,$
then $a^{\ast }=0,$ so $a=a^{\ast \ast }=0^{\ast }=1.\Box \medskip $

From the above, we deduce the following result.\medskip

\textbf{Corollary 35.}

\textit{1) A finite BL-algebra} $L$ \textit{with} $\max \mathcal{D}(L)\neq
0,1$ \textit{is not an MV-algebra}$.$

\textit{2) A finite MV-algebra} $L$ \textit{is not a chain iff} $\mathcal{D}%
(L)=\{0\};$

\textit{3) An finite MV-algebra that is not a chain is not a comet.}$\Box
\medskip $

\textbf{Proposition 36. ([\textbf{CFDP; 22}]) } \textit{If }$A$\textit{\ is
a finite commutative ring with }$\left\vert A\right\vert =n=p_{1}^{\alpha
_{1}}\cdot ...\cdot p_{r}^{\alpha _{r}},$\textit{\ then its set of ideals is
an MV-algebra. Of all its representations, only if }$A$\textit{\ is
isomorphic to the ring }$\underset{\alpha _{1}-time}{\underbrace{\mathbb{Z}%
_{p_{1}}\times \mathbb{Z}_{p_{1}}\times ...\times \mathbb{Z}_{p_{1}}}\times
...\times }\underset{\alpha _{r}-time}{\underbrace{\mathbb{Z}_{p_{r}}\times 
\mathbb{Z}_{p_{r}}\times ...\times \mathbb{Z}_{p_{r}}}}$\textit{\ the
lattice of its ideals is a Boolean algebra.\medskip }

\textbf{Examples 37.}

1) To generate a BL-comet with $k+4$ elements, $k\geq 1,$ organized as a
lattice as in Figure 1,

\begin{center}
\begin{tikzpicture}[scale=1]

  \coordinate (Top) at (0,4);        
  \coordinate (A) at (-1.5,3);       
  \coordinate (B) at (1.5,3);        
  \coordinate (Ak) at (0,2);         

  \coordinate (Akminus1) at (0,1.4); 
  \coordinate (Dots1) at (0,0.9);    
  \coordinate (A1) at (0,0.4);       
  \coordinate (A0) at (0,-0.2);      

  \draw (Top) -- (A) -- (Ak) -- (B) -- (Top);

  \draw (Ak) -- (Akminus1);
  \draw (A1) -- (A0);

  \foreach \pt in {Top,A,B,Ak,Akminus1,A1,A0}
    \fill (\pt) circle (1pt);

  \node[above] at (Top) {1};
  \node[left] at (A) {a};
  \node[right] at (B) {b};
  \node[right] at (Ak) {$a_k$};
  \node[right] at (Akminus1) {$a_{k-1}$};
  \node at (Dots1) {$\vdots$};
  \node[right] at (A1) {$a_1$};
  \node[below] at (A0) {$a_0=0$};

\end{tikzpicture}
\end{center}

\begin{equation*}
Figure\text{ 1.}
\end{equation*}%
we consider the commutative rings $\ (\mathbb{Z}_{2^{k}},+,\cdot )$ and $(%
\mathbb{Z}_{2}\times \mathbb{Z}_{2},+,\cdot ).$

We recall that $(Id(Z_{2^{k}}),\cap ,+,\otimes ,\rightarrow
,0=\{0\},1=Z_{2^{k}})$ is the only MV-chain (up to an isomorphism) with $k+1$
elements, see [\textbf{CFDP; 22}].

The ring $(Z_{2^{k}},+,\cdot )$ has $k+1$ ideals: $I_{0}=\{0\},$ $I_{1}=%
\widehat{2^{k-1}}Z_{2^{k}},$ ..., $I_{k-2}=\widehat{2^{2}}Z_{2^{k}},$ $%
I_{k-1}=\widehat{2}Z_{p^{k}},$ $I_{k}=Z_{2^{k}}$ and $I_{0}\subseteq
I_{1}\subseteq I_{2}\subseteq ...\subseteq I_{k}.$

For every $i,j\in \{0,...,k\}$ we have%
\begin{equation*}
I_{i}\rightarrow I_{j}=Z_{2^{k}}\text{ if }i\leq j\text{ and }I_{k-i+j}\text{
otherwise}
\end{equation*}%
and 
\begin{equation*}
I_{i}\oplus I_{j}=Z_{2^{k}}\text{ if }k\leq i+j\text{ and }I_{i+j}\text{
otherwise.}
\end{equation*}

\bigskip Also, $I_{i}^{\ast }=Ann(I_{i})=$ $I_{k-i}$ for every $i\in
\{0,...,k\}.$

We deduce that $I_{i}\otimes I_{j}=(I_{i}^{\ast }\oplus I_{j}^{\ast })^{\ast
}=Ann(I_{k-i}\oplus I_{k-j})=Ann(Z_{2^{k}})$ if $k\leq (k-i)+(k-j)$ and $%
Ann(I_{(k-i)+(k-j)})$ otherwise.

We conclude that 
\begin{equation*}
I_{i}\otimes I_{j}=I_{0\text{ }}\text{ if }i+j\leq k\text{ and }I_{i+j-k}%
\text{ otherwise.}
\end{equation*}

For the ring $(\mathbb{Z}_{2}\times \mathbb{Z}_{2},+,\cdot )$ the lattice of
ideals is $Id\left( \mathbb{Z}_{2}\times \mathbb{Z}_{2}\right) =\{\left( 
\widehat{0},\widehat{0}\right) ,\{\left( \widehat{0},\widehat{0}\right)
,\left( \widehat{0},\widehat{1}\right) \},\{\left( \widehat{0},\widehat{0}%
\right) ,\left( \widehat{1},\widehat{0}\right) \},\mathbb{Z}_{2}\times 
\mathbb{Z}_{2}\}=\{O,A,B,E\}$, which is a Boolean algebra $(Id\left( \mathbb{%
Z}_{2}\times \mathbb{Z}_{2}\right) ,\cap ,+,\otimes \rightarrow ,0=\{\left( 
\widehat{0},\widehat{0}\right) \},1=\mathbb{Z}_{2}\times \mathbb{Z}_{2}),$
so a BL-algebra, with the following operations: 
\begin{equation*}
\begin{tabular}{l|llll}
$\rightarrow $ & $O$ & $A$ & $B$ & $E$ \\ \hline
$O$ & $E$ & $E$ & $E$ & $E$ \\ 
$A$ & $B$ & $E$ & $B$ & $E$ \\ 
$B$ & $A$ & $A$ & $E$ & $E$ \\ 
$E$ & $O$ & $A$ & $B$ & $E$%
\end{tabular}%
\text{ and }%
\begin{tabular}{l|llll}
$\otimes $ & $O$ & $A$ & $B$ & $E$ \\ \hline
$O$ & $O$ & $O$ & $O$ & $O$ \\ 
$A$ & $O$ & $A$ & $O$ & $A$ \\ 
$B$ & $O$ & $O$ & $B$ & $B$ \\ 
$E$ & $O$ & $A$ & $B$ & $E$%
\end{tabular}%
\text{ }.
\end{equation*}%
If we consider two BL-algebras isomorphic with $(Id\left( \mathbb{Z}%
_{2^{k}}\right) ,\cap ,+,\otimes \rightarrow ,0=\{0\},1=\mathbb{Z}_{2^{k}})$
and $(Id\left( \mathbb{Z}_{2}\times \mathbb{Z}_{2}\right) ,\cap ,+,\otimes
\rightarrow ,0=\{\left( \widehat{0},\widehat{0}\right) \},1=\mathbb{Z}%
_{2}\times \mathbb{Z}_{2}),$ denoted by $\mathcal{L}_{1}=(L_{1}=%
\{0=a_{0},a_{1,}...a_{k}\},\wedge _{1},\vee _{1},\odot _{1},\rightarrow
_{1},0,a_{k})$ and $\mathcal{L}_{2}=(L_{2}=\{a_{k},a,b,1\},\wedge _{2},\vee
_{2},\odot _{2},\rightarrow _{2},a_{k},1),$ we can generate a BL-comet $%
\mathcal{L}_{1}\boxplus \mathcal{L}_{2}=(L_{1}\cup
L_{2}=\{0=a_{0},a_{1,}...a_{k},,a,b,1\},\wedge ,\vee ,\odot ,\rightarrow
,0,1)$ with $k+4$ elements, for any $k\geq 1.$

For example, for $k=4$ \ we obtain a BL-comet $\mathcal{L}_{1}\boxplus 
\mathcal{L}_{2}=(\{0=a_{0},a_{1,}a_{2,}a_{3},a_{4},a,b,1\},\wedge ,\vee
,\odot ,\rightarrow ,0,1)$ with the following operations:

\begin{equation*}
\text{ }%
\begin{tabular}{l|llllllll}
$\rightarrow $ & $0$ & $a_{1}$ & $a_{2}$ & $a_{3}$ & $a_{4}$ & $a$ & $b$ & $%
1 $ \\ \hline
$0$ & $1$ & $1$ & $1$ & $1$ & $1$ & $1$ & $1$ & $1$ \\ 
$a_{1}$ & $a_{3}$ & $1$ & $1$ & $1$ & $1$ & $1$ & $1$ & $1$ \\ 
$a_{2}$ & $a_{2}$ & $a_{3}$ & $1$ & $1$ & $1$ & $1$ & $1$ & $1$ \\ 
$a_{3}$ & $a_{1}$ & $a_{2}$ & $a_{3}$ & $1$ & $1$ & $1$ & $1$ & $1$ \\ 
$a_{4}$ & $0$ & $a_{1}$ & $a_{2}$ & $a_{3}$ & $1$ & $1$ & $1$ & $1$ \\ 
$a$ & $0$ & $a_{1}$ & $a_{2}$ & $a_{3}$ & $b$ & $1$ & $b$ & $1$ \\ 
$b$ & $0$ & $a_{1}$ & $a_{2}$ & $a_{3}$ & $a$ & $a$ & $1$ & $1$ \\ 
$1$ & $0$ & $a_{1}$ & $a_{2}$ & $a_{3}$ & $a_{4}$ & $a$ & $b$ & $1$%
\end{tabular}%
\text{ and }%
\begin{tabular}{l|llllllll}
$\odot $ & $0$ & $a_{1}$ & $a_{2}$ & $a_{3}$ & $a_{4}$ & $a$ & $b$ & $1$ \\ 
\hline
$0$ & $0$ & $0$ & $0$ & $0$ & $0$ & $0$ & $0$ & $0$ \\ 
$a_{1}$ & $0$ & $0$ & $0$ & $0$ & $a_{1}$ & $a_{1}$ & $a_{1}$ & $a_{1}$ \\ 
$a_{2}$ & $0$ & $0$ & $0$ & $a_{1}$ & $a_{2}$ & $a_{2}$ & $a_{2}$ & $a_{2}$
\\ 
$a_{3}$ & $0$ & $0$ & $a_{1}$ & $a_{2}$ & $a_{3}$ & $a_{3}$ & $a_{3}$ & $%
a_{3}$ \\ 
$a_{4}$ & $0$ & $a_{1}$ & $a_{2}$ & $a_{3}$ & $a_{4}$ & $a_{4}$ & $a_{4}$ & $%
a_{4}$ \\ 
$a$ & $0$ & $a_{1}$ & $a_{2}$ & $a_{3}$ & $a_{4}$ & $a$ & $a_{4}$ & $a$ \\ 
$b$ & $0$ & $a_{1}$ & $a_{2}$ & $a_{3}$ & $a_{4}$ & $a_{4}$ & $b$ & $b$ \\ 
$1$ & $0$ & $a_{1}$ & $a_{2}$ & $a_{3}$ & $a_{4}$ & $a$ & $b$ & $1$%
\end{tabular}%
.
\end{equation*}

2) To generate a BL-comet with $k+6$ elements, $k\geq 1$, organized as a
lattice as in Figure 2,

\begin{center}
\begin{tikzpicture}[scale=1.2]

  \coordinate (Top) at (0,4);        
  \coordinate (B) at (-1.5,3);       
  \coordinate (D) at (1.5,3);        
  \coordinate (A) at (0,2);          
  \coordinate (C) at (3,2);          
  \coordinate (Ak) at (1.5,1);       

  \coordinate (Akminus1) at (1.5,0.6); 
  \coordinate (Dots1) at (1.5,0.3);    
  \coordinate (A1) at (1.5,0);         
  \coordinate (A0) at (1.5,-0.4);      

  \draw (Top) -- (B) -- (A) -- (D) -- (Top);  
  \draw (D) -- (C) -- (Ak) -- (A) -- (D);     

  \draw (Ak) -- (Akminus1);
  \draw (A1) -- (A0);

  \foreach \pt in {Top,B,D,A,C,Ak,Akminus1,A1,A0}
    \fill (\pt) circle (1pt);

  \node[above] at (Top) {1};
  \node[left] at (B) {b};
  \node[right] at (D) {d};
  \node[left] at (A) {a};
  \node[right] at (C) {c};
  \node[right] at (Ak) {$a_k$};
  \node[right] at (Akminus1) {$a_{k-1}$};
  \node at (Dots1) {$\vdots$};
  \node[right] at (A1) {$a_1$};
  \node[below] at (A0) {$a_0=0$};

\end{tikzpicture}
\end{center}

\begin{equation*}
Figure\text{ 2.}
\end{equation*}

we consider the commutative rings $\ (\mathbb{Z}_{2^{k}},+,\cdot )$ and $(%
\mathbb{Z}_{2}\times \mathbb{Z}_{4},+,\cdot ).$

The ring $(Z_{2^{k}},+,\cdot )$ has $k+1$ ideals and $(Id\left( \mathbb{Z}%
_{2^{k}}\right) ,\cap ,+,\otimes \rightarrow ,0=\{0\},1=\mathbb{Z}_{2^{k}})$
is a BL-chain.

For $\mathbb{Z}_{2}\times \mathbb{Z}_{4}=\{\left( \widehat{0},\overline{0}%
\right) $, $\left( \widehat{0},\overline{1}\right) $, $\left( \widehat{0},%
\overline{2}\right) $, $\left( \widehat{0},\overline{3}\right) $, $\left( 
\widehat{1},\overline{0}\right) $, $\left( \widehat{1},\overline{1}\right) $%
, $\left( \widehat{1},\overline{2}\right) $, $\left( \widehat{1},\overline{3}%
\right) \}$, the lattice of ideals is \newline
$Id\left( \mathbb{Z}_{2}\times \mathbb{Z}_{4}\right) $=$\{\left( \widehat{0},%
\overline{0}\right) $,$\{\left( \widehat{0},\overline{0}\right) $,$\left( 
\widehat{0},\overline{1}\right) $,$\left( \widehat{0},\overline{2}\right) $,$%
\left( \widehat{0},\overline{3}\right) \}$,\newline
$\{\left( \widehat{0},\overline{0}\right) ,\left( \widehat{1},\overline{0}%
\right) ,\left( \widehat{0},\overline{2}\right) ,\left( \widehat{1},%
\overline{2}\right) \}$,$\{\left( \widehat{0},\overline{0}\right) ,\left( 
\widehat{0},\overline{2}\right) \}$, $\{\left( \widehat{0},\overline{0}%
\right) ,\left( \widehat{1},\overline{0}\right) \}$, $\mathbb{Z}_{2}\times 
\mathbb{Z}_{4}\}=\{O,B,D,A,C,E\}$ is an MV-algebra, with the following
operations:%
\begin{equation*}
\begin{tabular}{l|llllll}
$\rightarrow $ & $O$ & $A$ & $B$ & $C$ & $D$ & $E$ \\ \hline
$O$ & $E$ & $E$ & $E$ & $E$ & $E$ & $E$ \\ 
$A$ & $D$ & $E$ & $E$ & $D$ & $E$ & $E$ \\ 
$B$ & $C$ & $D$ & $E$ & $C$ & $D$ & $E$ \\ 
$C$ & $B$ & $B$ & $B$ & $E$ & $E$ & $E$ \\ 
$D$ & $A$ & $B$ & $B$ & $D$ & $E$ & $E$ \\ 
$E$ & $O$ & $A$ & $B$ & $C$ & $D$ & $E$%
\end{tabular}%
\text{ and }%
\begin{tabular}{l|llllll}
$\otimes $ & $O$ & $A$ & $B$ & $C$ & $D$ & $E$ \\ \hline
$O$ & $O$ & $O$ & $O$ & $O$ & $O$ & $O$ \\ 
$A$ & $O$ & $O$ & $A$ & $O$ & $O$ & $A$ \\ 
$B$ & $O$ & $A$ & $B$ & $O$ & $A$ & $B$ \\ 
$C$ & $O$ & $O$ & $O$ & $C$ & $C$ & $C$ \\ 
$D$ & $O$ & $O$ & $A$ & $C$ & $C$ & $D$ \\ 
$E$ & $O$ & $A$ & $B$ & $C$ & $D$ & $E$%
\end{tabular}%
\end{equation*}%
If we consider two BL-algebras isomorphic with $(Id\left( \mathbb{Z}%
_{2^{k}}\right) ,\cap ,+,\otimes \rightarrow ,0=\{0\},1=\mathbb{Z}_{2^{k}})$
and $(Id\left( \mathbb{Z}_{2}\times \mathbb{Z}_{4}\right) ,\cap ,+,\otimes
\rightarrow ,0=\{\left( \widehat{0},\overline{0}\right) \},1=\mathbb{Z}%
_{2}\times \mathbb{Z}_{4}),$ denoted by $\mathcal{L}_{1}=(L_{1}=%
\{0=a_{0},a_{1,}...a_{k}\},\wedge _{1},\vee _{1},\odot _{1},\rightarrow
_{1},0,a_{k})$ and $\mathcal{L}_{2}=(L_{2}=\{a_{k},a,b,c,d,1\},\wedge
_{2},\vee _{2},\odot _{2},\rightarrow _{2},a_{k},1),$ we can generate a
BL-comet $\mathcal{L}_{1}\boxplus \mathcal{L}_{2}=(L_{1}\cup
L_{2}=\{0=a_{0},a_{1,}...a_{k},a,b,c,d,1\},\wedge ,\vee ,\odot ,\rightarrow
,0,1)$ with $k+6$ elements, for any $k\geq 1.$

3) To generate a BL-comet with $k+8$ elements, $k\geq 1$, organized as a
lattice as in Figure 3,

\begin{center}
\begin{tikzpicture}[scale=1.3]

  \coordinate (Top) at (0,3);       
  \coordinate (U) at (-1.1,2.2);    
  \coordinate (V) at (1.1,2.2);     
  \coordinate (Z) at (-0.1,1.5);    
  \coordinate (T) at (0.1,1.5);     
  \coordinate (X) at (-1.1,0.9);    
  \coordinate (Y) at (1.1,0.9);     
  \coordinate (Ak) at (0.1,0.3);      

  \coordinate (Akminus1) at (0.1,-0.2);    
  \coordinate (Dots1) at (0.1,-0.6);       
  \coordinate (A1) at (0.1,-1.0);          
  \coordinate (A0) at (0.1,-1.4);          

  \draw (Top)--(U)--(X)--(Ak)--(Y)--(V)--(Top);
  \draw (Z)--(Top);
  \draw (Z)--(X);
  \draw (T)--(U);
  \draw (T)--(V);
  \draw (Z)--(Y);
  \draw (T)--(Ak);

  \draw (Ak)--(Akminus1);
  \draw (A1)--(A0);

  \foreach \pt in {Top,U,V,Z,T,X,Y,Ak,Akminus1,A1,A0, T}
    \fill (\pt) circle (1pt);

  \node[above] at (Top) {1};
  \node[left] at (U) {z};
  \node[right] at (V) {v};
  \node[left] at (Z) {u};
  \node[right] at (T) {b};
  \node[left] at (X) {x};
  \node[right] at (Y) {t};
  \node[right] at (Ak) {$a_k$};
  \node[right] at (Akminus1) {$a_{k-1}$};
  \node at (Dots1) {$\vdots$};
  \node[right] at (A1) {$a_1$};
  \node[below] at (A0) {$a_0=0$};

\end{tikzpicture}
\end{center}

\begin{equation*}
Figure\text{ 3.}
\end{equation*}%
we consider the commutative rings $\ (\mathbb{Z}_{2^{k}},+,\cdot )$ and $(%
\mathbb{Z}_{2}\times \mathbb{Z}_{2}\times \mathbb{Z}_{2},+,\cdot ).$

The ring $(Z_{2^{k}},+,\cdot )$ has $k+1$ ideals and $(Id\left( \mathbb{Z}%
_{2^{k}}\right) ,\cap ,+,\otimes \rightarrow ,0=\{0\},1=\mathbb{Z}_{2^{k}})$
is a BL-chain.

For $\mathbb{Z}_{2}\times \mathbb{Z}_{2}\times \mathbb{Z}_{2}$ the lattice
of ideals $Id\left( \mathbb{Z}_{2}\times \mathbb{Z}_{2}\times \mathbb{Z}%
_{2}\right) $ has $8$ ideals denoted $\{O,X,Y,Z,T,U,V,E\}$ and is a Boolean
algebra with the following operations:

\begin{equation*}
\text{ }%
\begin{tabular}{l|llllllll}
$\rightarrow $ & $O$ & $X$ & $Y$ & $Z$ & $T$ & $U$ & $V$ & $E$ \\ \hline
$O$ & $E$ & $E$ & $E$ & $E$ & $E$ & $E$ & $E$ & $E$ \\ 
$X$ & $V$ & $E$ & $V$ & $E$ & $V$ & $E$ & $V$ & $E$ \\ 
$Y$ & $U$ & $U$ & $E$ & $E$ & $U$ & $U$ & $E$ & $E$ \\ 
$Z$ & $T$ & $U$ & $V$ & $E$ & $T$ & $U$ & $V$ & $E$ \\ 
$T$ & $Z$ & $Z$ & $Z$ & $Z$ & $E$ & $E$ & $E$ & $E$ \\ 
$U$ & $Y$ & $Z$ & $Y$ & $Z$ & $V$ & $E$ & $V$ & $E$ \\ 
$V$ & $X$ & $X$ & $Z$ & $Z$ & $U$ & $U$ & $E$ & $E$ \\ 
$E$ & $O$ & $X$ & $Y$ & $Z$ & $T$ & $U$ & $V$ & $E$%
\end{tabular}%
\text{ and }%
\begin{tabular}{l|llllllll}
$\otimes $ & $O$ & $X$ & $Y$ & $Z$ & $T$ & $U$ & $V$ & $E$ \\ \hline
$O$ & $O$ & $O$ & $O$ & $O$ & $O$ & $O$ & $O$ & $O$ \\ 
$X$ & $O$ & $X$ & $O$ & $X$ & $O$ & $X$ & $O$ & $X$ \\ 
$Y$ & $O$ & $O$ & $Y$ & $Y$ & $O$ & $O$ & $Y$ & $Y$ \\ 
$Z$ & $O$ & $X$ & $Y$ & $Z$ & $O$ & $X$ & $Y$ & $Z$ \\ 
$T$ & $O$ & $O$ & $O$ & $O$ & $T$ & $T$ & $T$ & $T$ \\ 
$U$ & $O$ & $X$ & $O$ & $X$ & $T$ & $U$ & $T$ & $U$ \\ 
$V$ & $O$ & $O$ & $Y$ & $Y$ & $T$ & $T$ & $V$ & $V$ \\ 
$E$ & $O$ & $X$ & $Y$ & $Z$ & $T$ & $U$ & $V$ & $E$%
\end{tabular}%
\end{equation*}

If we consider two BL-algebras isomorphic with $(Id\left( \mathbb{Z}%
_{2^{k}}\right) ,\cap ,+,\otimes \rightarrow ,0=\{0\},1=\mathbb{Z}_{2^{k}})$
and $(Id\left( \mathbb{Z}_{2}\times \mathbb{Z}_{2}\times \mathbb{Z}%
_{2}\right) ,\cap ,+,\otimes \rightarrow ,0=\{\left( \widehat{0},\widehat{0}%
\right) \},1=\mathbb{Z}_{2}\times \mathbb{Z}_{2}\times \mathbb{Z}_{2}),$
denoted by $\mathcal{L}_{1}=(L_{1}=\{0=a_{0},a_{1,}...a_{k}\},\wedge
_{1},\vee _{1},\odot _{1},\rightarrow _{1},0,a_{k})$ and $\mathcal{L}%
_{2}=(L_{2}=\{a_{k},x,y,z,t,u,v,1\},\wedge _{2},\vee _{2},\odot
_{2},\rightarrow _{2},a_{k},1),$ we can generate a BL-comet $\mathcal{L}%
_{1}\boxplus \mathcal{L}_{2}=(L_{1}\cup
L_{2}=\{0=a_{0},a_{1,}...a_{k},x,y,z,t,u,v,1\},\wedge ,\vee ,\odot
,\rightarrow ,0,1)$ with $k+8$ elements, for any $k\geq 1.$

4) To generate a BL-comet with $k+9$ elements, $k\geq 1$, organized as a
lattice as in Figure 4,

\begin{center}
\begin{tikzpicture}[scale=1.3]

  \coordinate (Top) at (0,4);        
  \coordinate (Y) at (-1.5,3);       
  \coordinate (V) at (1.5,3);        
  \coordinate (Ak) at (0,2);         
  \coordinate (T) at (0,3); 
  \coordinate (U) at (-0.7,3.5); 
  \coordinate (W) at (0.7,3.5); 
  \coordinate (X) at (-0.7,2.5); 
  \coordinate (Z) at (0.7,2.5); 

  \coordinate (Akminus1) at (0,1.4); 
  \coordinate (Dots1) at (0,0.9);    
  \coordinate (A1) at (0,0.4);       
  \coordinate (A0) at (0,-0.2);      

  \draw (Top) -- (Y) -- (Ak) -- (V) -- (Top);

  \draw (Ak) -- (Akminus1);
  \draw (A1) -- (A0);
  \draw (U) -- (T);
  \draw (X) -- (T);
  \draw (W) -- (T);
  \draw (Z) -- (T);

  \foreach \pt in {Top,Y,V,Ak,Akminus1,A1,A0,T,U,W,X,Z}
    \fill (\pt) circle (1pt);

  \node[above] at (Top) {1};
  \node[left] at (Y) {y};
  \node[right] at (V) {v};
  \node[right] at (Ak) {$a_k$};
  \node[right] at (Akminus1) {$a_{k-1}$};
  \node at (Dots1) {$\vdots$};
  \node[right] at (A1) {$a_1$};
  \node[below] at (A0) {$a_0=0$};
  \node[right] at  (T) {t};
  \node[left] at (U) {u};
  \node[right] at (W) {w};
  \node[left] at (X) {x};
  \node[right] at (Z) {z};

\end{tikzpicture}
\end{center}

\begin{equation*}
Figure\text{ 4.}
\end{equation*}%
we consider the commutative rings $\ (\mathbb{Z}_{2^{k}},+,\cdot )$ and $(%
\mathbb{Z}_{4}\times \mathbb{Z}_{4},+,\cdot ).$

$Id(Z_{2^{k}})$ is a BL-chain with $k+1$ elements and $Id\left( \mathbb{Z}%
_{4}\times \mathbb{Z}_{4}\right) $ is an MV-algebra with $9$ elements
denoted $\{O,X,Y,Z,T,U,V,W,E\}$ with the following operations:%
\begin{equation*}
\begin{tabular}{l|lllllllll}
$\rightarrow $ & $O$ & $X$ & $Y$ & $Z$ & $T$ & $U$ & $V$ & $W$ & $E$ \\ 
\hline
$O$ & $E$ & $E$ & $E$ & $E$ & $E$ & $E$ & $E$ & $E$ & $E$ \\ 
$X$ & $W$ & $E$ & $E$ & $W$ & $E$ & $E$ & $W$ & $E$ & $E$ \\ 
$Y$ & $V$ & $W$ & $E$ & $V$ & $W$ & $E$ & $V$ & $W$ & $E$ \\ 
$Z$ & $U$ & $U$ & $U$ & $E$ & $E$ & $E$ & $E$ & $E$ & $E$ \\ 
$T$ & $T$ & $U$ & $U$ & $W$ & $E$ & $E$ & $W$ & $E$ & $E$ \\ 
$U$ & $Z$ & $T$ & $U$ & $V$ & $W$ & $E$ & $V$ & $W$ & $E$ \\ 
$V$ & $Y$ & $Y$ & $Y$ & $U$ & $U$ & $U$ & $E$ & $E$ & $E$ \\ 
$W$ & $X$ & $Y$ & $Y$ & $T$ & $U$ & $U$ & $W$ & $E$ & $E$ \\ 
$E$ & $O$ & $X$ & $Y$ & $Z$ & $T$ & $U$ & $V$ & $W$ & $E$%
\end{tabular}%
\text{ and }
\end{equation*}%
\begin{equation*}
\begin{tabular}{l|lllllllll}
$\otimes $ & $O$ & $X$ & $Y$ & $Z$ & $T$ & $U$ & $V$ & $W$ & $E$ \\ \hline
$O$ & $O$ & $O$ & $O$ & $O$ & $O$ & $O$ & $O$ & $O$ & $O$ \\ 
$X$ & $O$ & $O$ & $X$ & $O$ & $O$ & $X$ & $O$ & $O$ & $X$ \\ 
$Y$ & $O$ & $X$ & $Y$ & $O$ & $X$ & $Y$ & $O$ & $X$ & $Y$ \\ 
$Z$ & $O$ & $O$ & $O$ & $O$ & $O$ & $O$ & $Z$ & $Z$ & $Z$ \\ 
$T$ & $O$ & $O$ & $X$ & $O$ & $O$ & $X$ & $Z$ & $Z$ & $T$ \\ 
$U$ & $O$ & $X$ & $Y$ & $O$ & $X$ & $Y$ & $Z$ & $T$ & $U$ \\ 
$V$ & $O$ & $O$ & $O$ & $Z$ & $Z$ & $Z$ & $V$ & $V$ & $V$ \\ 
$W$ & $O$ & $O$ & $X$ & $Z$ & $Z$ & $T$ & $V$ & $V$ & $W$ \\ 
$E$ & $O$ & $X$ & $Y$ & $Z$ & $T$ & $U$ & $V$ & $W$ & $E$%
\end{tabular}%
\end{equation*}

If we consider two BL-algebras isomorphic with $(Id\left( \mathbb{Z}%
_{2^{k}}\right) ,\cap ,+,\otimes \rightarrow ,0=\{0\},1=\mathbb{Z}_{2^{k}})$
and $(Id\left( \mathbb{Z}_{4}\times \mathbb{Z}_{4}\right) ,\cap ,+,\otimes
\rightarrow ,0=\{\left( \widehat{0},\widehat{0}\right) \},1=\mathbb{Z}%
_{4}\times \mathbb{Z}_{4}),$ denoted by $\mathcal{L}_{1}=(L_{1}=%
\{0=a_{0},a_{1,}...a_{k}\},\wedge _{1},\vee _{1},\odot _{1},\rightarrow
_{1},0,a_{k})$ and $\mathcal{L}_{2}=(L_{2}=\{a_{k},x,y,z,t,u,v,w,1\},\wedge
_{2},\vee _{2},\odot _{2},\rightarrow _{2},a_{k},1),$ we can generate a
BL-comet $\mathcal{L}_{1}\boxplus \mathcal{L}_{2}=(L_{1}\cup
L_{2}=\{0=a_{0},a_{1,}...a_{k},x,y,z,t,u,v,w,1\},\wedge ,\vee ,\odot
,\rightarrow ,0,1)$ with $k+9$ elements, for any $k\geq 1.\medskip $

\textbf{Remark 38.}Using Example 37, for any $n\geq 5$, we can generate
BL-comets with $n$ elements which are not chains.

In [BV;10], isomorphism classes of BL-algebras of size $n\leq 12$ were just
counted, not constructed, using computer algorithms. Up to an isomorphism,
there are $1$ BL-algebra of size $2$, $2$ BL-algebras of size $3$, $5$
BL-algebras of size $4$, $9$ BL-algebras of size $5$, $20$ BL-algebras of
size $6$, $38$ BL-algebras of size $7$, $81$ BL-algebras of size $8$, $160$
BL-algebras of size $9$, $326$ BL-algebras of size $10$, $643$ BL-algebra of
size $11$ and $1314$ BL-algebras of size $12$. In [FP; 22] we construct (up
to an isomorphism) all finite BL-algebras with $2\leq n\leq 5$
elements.\medskip

\textbf{Table 1} present a summary of the structure of BL-algebras $L$ with $%
2\leq n\leq 5$ elements:

\begin{equation*}
\text{\textbf{Table 1:}}
\end{equation*}

\textbf{\ }%
\begin{tabular}{lll}
$\left\vert L\right\vert \mathbf{=n}$ & \textbf{Nr of BL-alg } & \textbf{%
Structure } \\ 
$n=2$ & $1$ & $\left\{ Id(\mathbb{Z}_{2})\text{ (chain, MV, COMET)}\right. $
\\ 
$n=3$ & $2$ & $\left\{ 
\begin{array}{c}
Id(\mathbb{Z}_{4})\text{ (chain, MV, COMET)} \\ 
Id(\mathbb{Z}_{2})\boxplus Id(\mathbb{Z}_{2})\text{ (chain, BL, COMET)}%
\end{array}%
\right. $ \\ 
$n=4$ & $5$ & $\left\{ 
\begin{array}{c}
Id(\mathbb{Z}_{8})\text{ (chain, MV, COMET)} \\ 
Id(\mathbb{Z}_{2}\times \mathbb{Z}_{2})\text{ (MV, NOT COMET)} \\ 
Id(\mathbb{Z}_{2})\boxplus Id(\mathbb{Z}_{4})\text{ (chain, BL, COMET)} \\ 
Id(\mathbb{Z}_{4})\boxplus Id(\mathbb{Z}_{2})\text{ \ (chain, BL, COMET)} \\ 
Id(\mathbb{Z}_{2})\boxplus (Id(\mathbb{Z}_{2})\boxplus Id(\mathbb{Z}_{2}))%
\text{ (chain, BL, COMET)}%
\end{array}%
\right. $ \\ 
$n=5$ & $9$ & $\left\{ 
\begin{array}{c}
Id(\mathbb{Z}_{16})\text{ (chain, MV, COMET)} \\ 
Id(\mathbb{Z}_{2})\boxplus Id(\mathbb{Z}_{8})\text{ (chain, BL, COMET)} \\ 
Id(\mathbb{Z}_{2})\boxplus Id(\mathbb{Z}_{2}\times \mathbb{Z}_{2})\text{
(BL, COMET)} \\ 
Id(\mathbb{Z}_{2})\boxplus (Id(\mathbb{Z}_{2})\boxplus Id(\mathbb{Z}_{4}))%
\text{ (chain, BL, COMET)} \\ 
Id(\mathbb{Z}_{2})\boxplus (Id(\mathbb{Z}_{4})\boxplus Id(\mathbb{Z}_{2}))%
\text{ (chain, BL, COMET)} \\ 
Id(\mathbb{Z}_{2})\boxplus (Id(\mathbb{Z}_{2})\boxplus (Id(\mathbb{Z}%
_{2})\boxplus Id(\mathbb{Z}_{2})))\text{ (chain, BL, COMET)} \\ 
Id(\mathbb{Z}_{8})\boxplus Id(\mathbb{Z}_{2})\text{ (chain, BL, COMET)} \\ 
(Id(\mathbb{Z}_{4})\boxplus Id(\mathbb{Z}_{2}))\boxplus Id(\mathbb{Z}_{2})%
\text{ (chain, BL, COMET)} \\ 
Id(\mathbb{Z}_{4})\boxplus Id(\mathbb{Z}_{4})\text{ (chain, BL, COMET)}%
\end{array}%
\right. $%
\end{tabular}

In the following, by using the ordinal sum of two BL-algebras we generate
all (up to an isomorphism) finite BL-algebras (which are not MV-algebras )
with $n=6$ elements. This method can be used to construct finite BL-algebras
of larger size, the inconvenience being the large number of BL-algebras that
should be generated.\medskip

\textbf{Theorem 39.} \textit{i)} \textit{All BL-algebras with} $6$ \textit{%
elements (which are not MV-algebras) \ can be generated as ordinal sum }$%
\mathcal{L}_{1}\boxplus \mathcal{L}_{2}$ \textit{of }$\ $\textit{two
BL-algebras} $\mathcal{L}_{1}$ \textit{and} $\mathcal{L}_{2}$ \textit{in the
following ways}:%
\begin{equation*}
\mathcal{L}_{1}\text{ \textit{is a BL-chain with} }2\text{ \textit{elements
and} }\mathcal{L}_{2}\text{ \textit{is a BL-algebra with} }5\text{ \textit{%
elements,} }
\end{equation*}%
\textit{or}%
\begin{equation*}
\mathcal{L}_{1}\text{ \textit{is a BL-chain with} }3\text{ \textit{elements
and} }\mathcal{L}_{2}\text{ \textit{is a BL-algebra with} }4\text{ \textit{%
elements, }}
\end{equation*}%
\textit{or}%
\begin{equation*}
\mathcal{L}_{1}\text{ \textit{is a BL-chain with} }4\text{ \textit{elements
and} }\mathcal{L}_{2}\text{ \textit{is a BL-algebra with} }3\text{ \textit{%
elements,} }
\end{equation*}%
\textit{or}%
\begin{equation*}
\mathcal{L}_{1}\text{ \textit{is a BL-chain with} }5\text{ \textit{elements
and} }\mathcal{L}_{2}\text{ \textit{is a BL-algebra with} }2\text{ \textit{%
elements.} }
\end{equation*}

\textit{ii) All }$18$ \textit{BL-algebras with} $6$ \textit{elements that
are not MV-algebras are BL-comets.}

\textit{iii) There are} $20$ \textit{BL-algebras with} $6$ \textit{%
elements.\medskip }

\textbf{Proof. i) Case 1.}%
\begin{equation*}
\mathcal{L}_{1}\text{ is a BL-chain with }2\text{ elements and }\mathcal{L}%
_{2}\text{ is a BL-algebra with }5\text{ elements. }
\end{equation*}%
We obtain the following BL-algebras: 
\begin{eqnarray*}
&&Id(\mathbb{Z}_{2})\boxplus Id(\mathbb{Z}_{16}),\text{ }Id(\mathbb{Z}%
_{2})\boxplus \lbrack Id(\mathbb{Z}_{2})\boxplus Id(\mathbb{Z}_{8})],Id(%
\mathbb{Z}_{2})\boxplus \lbrack Id(\mathbb{Z}_{2})\boxplus Id(\mathbb{Z}%
_{2}\times \mathbb{Z}_{2})], \\
&&Id(\mathbb{Z}_{2})\boxplus \lbrack Id(\mathbb{Z}_{2})\boxplus (Id(\mathbb{Z%
}_{2})\boxplus Id(\mathbb{Z}_{4}))],\text{ }Id(\mathbb{Z}_{2})\boxplus
\lbrack Id(\mathbb{Z}_{2})\boxplus (Id(\mathbb{Z}_{4})\boxplus Id(\mathbb{Z}%
_{2}))], \\
&&Id(\mathbb{Z}_{2})\boxplus \{Id(\mathbb{Z}_{2})\boxplus \lbrack Id(\mathbb{%
Z}_{2})\boxplus (Id(\mathbb{Z}_{2})\boxplus Id(\mathbb{Z}_{2}))]\},\text{ }%
Id(\mathbb{Z}_{2})\boxplus (Id(\mathbb{Z}_{8})\boxplus Id(\mathbb{Z}_{2})),
\\
&&Id(\mathbb{Z}_{2})\boxplus \lbrack (Id(\mathbb{Z}_{4})\boxplus Id(\mathbb{Z%
}_{2}))\boxplus Id(\mathbb{Z}_{2})],\text{ }Id(\mathbb{Z}_{2})\boxplus
\lbrack Id(\mathbb{Z}_{4})\boxplus Id(\mathbb{Z}_{4})].
\end{eqnarray*}%
\textbf{Case 2.}%
\begin{equation*}
\mathcal{L}_{1}\text{ is a BL-chain with }3\text{ elements and }\mathcal{L}%
_{2}\text{ is a BL-algebra with }4\text{ elements. }
\end{equation*}%
We obtain the following BL-algebras: 
\begin{eqnarray*}
&&Id(\mathbb{Z}_{4})\boxplus Id(\mathbb{Z}_{8}),\text{ }Id(\mathbb{Z}%
_{4})\boxplus Id(\mathbb{Z}_{2}\times \mathbb{Z}_{2}),\text{ } \\
&&Id(\mathbb{Z}_{4})\boxplus \lbrack Id(\mathbb{Z}_{2})\boxplus Id(\mathbb{Z}%
_{4})], \\
Id(\mathbb{Z}_{4})\boxplus \lbrack Id(\mathbb{Z}_{4})\boxplus Id(\mathbb{Z}%
_{2})], &&\text{ }Id(\mathbb{Z}_{4})\boxplus \lbrack Id(\mathbb{Z}%
_{2})\boxplus (Id(\mathbb{Z}_{2})\boxplus Id(\mathbb{Z}_{2}))],\text{ } \\
\lbrack Id(\mathbb{Z}_{2})\boxplus Id(\mathbb{Z}_{2})]\boxplus Id(\mathbb{Z}%
_{8}),\text{ }[ &&Id(\mathbb{Z}_{2})\boxplus Id(\mathbb{Z}_{2})]\boxplus Id(%
\mathbb{Z}_{2}\times \mathbb{Z}_{2}), \\
\text{ }[Id(\mathbb{Z}_{2})\boxplus Id(\mathbb{Z}_{2})]\boxplus \lbrack Id(%
\mathbb{Z}_{2})\boxplus Id(\mathbb{Z}_{4})],\text{ } &&[Id(\mathbb{Z}%
_{2})\boxplus Id(\mathbb{Z}_{2})]\boxplus \lbrack Id(\mathbb{Z}%
_{4}))\boxplus Id(\mathbb{Z}_{2})], \\
&&\text{ }[Id(\mathbb{Z}_{2})\boxplus Id(\mathbb{Z}_{2})]\boxplus \lbrack Id(%
\mathbb{Z}_{2})\boxplus (Id(\mathbb{Z}_{2})\boxplus Id(\mathbb{Z}_{2}))].
\end{eqnarray*}%
\textbf{Case 3.}%
\begin{equation*}
\mathcal{L}_{1}\text{ is a BL-chain with }4\text{ elements and }\mathcal{L}%
_{2}\text{ is a BL-algebra with }3\text{ elements. }
\end{equation*}%
We obtain the following BL-algebras: 
\begin{eqnarray*}
&&Id(\mathbb{Z}_{8})\boxplus Id(\mathbb{Z}_{4}),\text{ }Id(\mathbb{Z}%
_{8})\boxplus \lbrack Id(\mathbb{Z}_{2})\boxplus Id(\mathbb{Z}_{2})],[Id(%
\mathbb{Z}_{2})\boxplus Id(\mathbb{Z}_{4})]\boxplus Id(\mathbb{Z}_{4}), \\
&&[Id(\mathbb{Z}_{2})\boxplus Id(\mathbb{Z}_{4})]\boxplus \lbrack Id(\mathbb{%
Z}_{2})\boxplus Id(\mathbb{Z}_{2})],\text{ }[Id(\mathbb{Z}_{4})\boxplus Id(%
\mathbb{Z}_{2})]\boxplus Id(\mathbb{Z}_{4}), \\
&&[Id(\mathbb{Z}_{4})\boxplus Id(\mathbb{Z}_{2})]\boxplus \lbrack Id(\mathbb{%
Z}_{2})\boxplus Id(\mathbb{Z}_{2})],\text{ }[Id(\mathbb{Z}_{2})\boxplus (Id(%
\mathbb{Z}_{2})\boxplus Id(\mathbb{Z}_{2}))]\boxplus Id(\mathbb{Z}_{4}), \\
&&[Id(\mathbb{Z}_{2})\boxplus (Id(\mathbb{Z}_{2})\boxplus Id(\mathbb{Z}%
_{2}))]\boxplus \lbrack Id(\mathbb{Z}_{2})\boxplus Id(\mathbb{Z}_{2})].
\end{eqnarray*}%
\textbf{Case 4.}%
\begin{equation*}
\mathcal{L}_{1}\text{ is a BL-chain with }5\text{ elements and }\mathcal{L}%
_{2}\text{ is a BL-algebra with }2\text{ elements. }
\end{equation*}%
We obtain the following BL-algebras: 
\begin{eqnarray*}
&&Id(\mathbb{Z}_{16})\boxplus Id(\mathbb{Z}_{2}),\text{ }[Id(\mathbb{Z}%
_{2})\boxplus Id(\mathbb{Z}_{8})]\boxplus Id(\mathbb{Z}_{2}),[Id(\mathbb{Z}%
_{2})\boxplus (Id(\mathbb{Z}_{2})\boxplus Id(\mathbb{Z}_{4}))]\boxplus Id(%
\mathbb{Z}_{2}), \\
\lbrack &&Id(\mathbb{Z}_{2})\boxplus (Id(\mathbb{Z}_{4})\boxplus Id(\mathbb{Z%
}_{2}))]\boxplus Id(\mathbb{Z}_{2}),\text{ }[Id(\mathbb{Z}_{2})\boxplus (Id(%
\mathbb{Z}_{2})\boxplus (Id(\mathbb{Z}_{2})\boxplus Id(\mathbb{Z}%
_{2})))]\boxplus Id(\mathbb{Z}_{2}),\text{ } \\
&&[Id(\mathbb{Z}_{8})\boxplus Id(\mathbb{Z}_{2})]\boxplus Id(\mathbb{Z}%
_{2}),[(Id(\mathbb{Z}_{4})\boxplus Id(\mathbb{Z}_{2}))\boxplus Id(\mathbb{Z}%
_{2})]\boxplus Id(\mathbb{Z}_{2}),[Id(\mathbb{Z}_{4})\boxplus Id(\mathbb{Z}%
_{4})]\boxplus Id(\mathbb{Z}_{2}).
\end{eqnarray*}%
Since $\boxplus $ is associative, we obtain only $18$ BL-algebras of which $%
16$ are chains.

(ii). Obviously, see Table 2.

(iii). In addition, from all $18$ BL-algebras previously generated, there
are two MV-algebras: $Id(\mathbb{Z}_{32})$ and $Id(\mathbb{Z}_{2}\times 
\mathbb{Z}_{4}),$ see [\textbf{CFDP; 22}].$\Box \medskip $

\textbf{Table 2} present a summary of the structure of BL-algebras $L$ with $%
n=6$ elements:

\begin{equation*}
\text{\textbf{Table 2:}}
\end{equation*}

\begin{tabular}{ll}
$Id(\mathbb{Z}_{2})\boxplus Id(\mathbb{Z}_{16})$ & BL-chain $\Rightarrow $%
COMET \\ 
$Id(\mathbb{Z}_{2})\boxplus \lbrack Id(\mathbb{Z}_{2})\boxplus Id(\mathbb{Z}%
_{8})]$ & BL-chain $\Rightarrow $COMET \\ 
$Id(\mathbb{Z}_{2})\boxplus \lbrack Id(\mathbb{Z}_{2})\boxplus Id(\mathbb{Z}%
_{2}\times \mathbb{Z}_{2})]$ & BL $\Rightarrow $COMET, NOT CHAIN \\ 
$Id(\mathbb{Z}_{2})\boxplus \lbrack Id(\mathbb{Z}_{2})\boxplus (Id(\mathbb{Z}%
_{2})\boxplus Id(\mathbb{Z}_{4}))]$ & BL-chain $\Rightarrow $COMET \\ 
$\text{ }Id(\mathbb{Z}_{2})\boxplus \lbrack Id(\mathbb{Z}_{2})\boxplus (Id(%
\mathbb{Z}_{4})\boxplus Id(\mathbb{Z}_{2}))]$ & BL-chain $\Rightarrow $COMET
\\ 
$Id(\mathbb{Z}_{2})\boxplus \{Id(\mathbb{Z}_{2})\boxplus \lbrack Id(\mathbb{Z%
}_{2})\boxplus (Id(\mathbb{Z}_{2})\boxplus Id(\mathbb{Z}_{2}))]\}$ & 
BL-chain $\Rightarrow $COMET \\ 
$\text{ }Id(\mathbb{Z}_{2})\boxplus \lbrack Id(\mathbb{Z}_{8})\boxplus Id(%
\mathbb{Z}_{2})]$ & BL-chain $\Rightarrow $COMET \\ 
$Id(\mathbb{Z}_{2})\boxplus \lbrack (Id(\mathbb{Z}_{4})\boxplus Id(\mathbb{Z}%
_{2}))\boxplus Id(\mathbb{Z}_{2})]$ & BL-chain $\Rightarrow $COMET \\ 
$\text{ }Id(\mathbb{Z}_{2})\boxplus \lbrack Id(\mathbb{Z}_{4})\boxplus Id(%
\mathbb{Z}_{4})]$ & BL-chain $\Rightarrow $COMET \\ 
$Id(\mathbb{Z}_{4})\boxplus Id(\mathbb{Z}_{8})\text{ }$ & BL-chain $%
\Rightarrow $COMET \\ 
$Id(\mathbb{Z}_{4})\boxplus Id(\mathbb{Z}_{2}\times \mathbb{Z}_{2})$ & BL $%
\Rightarrow $COMET, NOT CHAIN \\ 
$Id(\mathbb{Z}_{4})\boxplus \lbrack Id(\mathbb{Z}_{2})\boxplus Id(\mathbb{Z}%
_{4})]$ & BL-chain $\Rightarrow $COMET \\ 
$Id(\mathbb{Z}_{4})\boxplus \lbrack Id(\mathbb{Z}_{4})\boxplus Id(\mathbb{Z}%
_{2})]\text{ }$ & BL-chain $\Rightarrow $COMET \\ 
$\text{ }Id(\mathbb{Z}_{4})\boxplus \lbrack Id(\mathbb{Z}_{2})\boxplus (Id(%
\mathbb{Z}_{2})\boxplus Id(\mathbb{Z}_{2}))]$ & BL-chain $\Rightarrow $COMET
\\ 
$Id(\mathbb{Z}_{8})\boxplus Id(\mathbb{Z}_{4})$ & BL-chain $\Rightarrow $%
COMET \\ 
$Id(\mathbb{Z}_{8})\boxplus \lbrack Id(\mathbb{Z}_{2})\boxplus Id(\mathbb{Z}%
_{4})]$ & BL-chain $\Rightarrow $COMET \\ 
$Id(\mathbb{Z}_{16})\boxplus Id(\mathbb{Z}_{2})$ & BL-chain $\Rightarrow $%
COMET \\ 
$\lbrack Id(\mathbb{Z}_{8})\boxplus Id(\mathbb{Z}_{2})]\boxplus Id(\mathbb{Z}%
_{2})$ & BL-chain $\Rightarrow $COMET \\ 
$Id(\mathbb{Z}_{2}\times \mathbb{Z}_{4})\text{ }$ & unordered MV $%
\Rightarrow $NOT COMET \\ 
$Id(\mathbb{Z}_{32})$ & MV-chain $\Rightarrow $COMET%
\end{tabular}%
.

\textbf{Corollary 40. }A finite BL-algebras with $n$ \textit{elements} ($%
n\leq 6)$ \textit{is not a comet iff it is an unordered MV-algebras.\medskip 
}

Finally, \textbf{Table 3 }present a summary for the number of MV-algebras,
MV-chains, BL-algebras, BL-chains and BL-comets with $n\leq 6$ elements:

\begin{equation*}
\text{\textbf{Table 3}}
\end{equation*}

\medskip

\begin{tabular}{llllll}
& $n=2$ & $n=3$ & $n=4$ & $n=5$ & $n=6$ \\ 
MV-algebras & 1 & 1 & 2 & 1 & 2 \\ 
MV-chains & 1 & 1 & 1 & 1 & 1 \\ 
BL-algebras & 1 & 2 & 5 & 9 & 20 \\ 
BL-chains & 1 & 2 & 4 & 8 & 17 \\ 
BL-comets & 1 & 2 & 3 & 9 & 19%
\end{tabular}

From the above results, we remark that a finite BL-algebra is a BL-comet or
an unordered MV-algebra, that means an MV-algebra which is not an MV-chain.
Now, we can state and demonstrate the main result of this paper.\medskip

\textbf{Theorem 41.} \textit{If} $L$ \textit{is a finite BL-algebra}, 
\textit{which is not an MV-algebra,} \textit{then there is no commutative
and unitary rings} $R$ \textit{such that} $Id\left( R\right) =L$.\medskip 

\textbf{Proof.} First, we prove the following Lemma.\medskip

\textbf{Lemma.} \textit{If} $R$ \textit{is a commutative, unitary and local
Artinian ring with a unique minimal ideal} $I_{m}$, \textit{then} $R$ 
\textit{is a chain ring}.\medskip

\textit{Proof of the Lemma.} Let \thinspace $R$ be a commuative and unitary
ring. The scole of the ring $R,Soc(R)$ is the sum of its minimal ideals. In
our case, $Soc\left( R\right) =I_{m}$. It is clear that $I_{m}$ is a
principal ideal, due to its minimality. We consider the ring $G=R/I_{m}$. We
have $Soc\left( G\right) =\sum \widehat{J},\widehat{J}$ minimal ideals in $%
R/I_{m}$. That means $J$ are those minimal ideals in $R$ containing $I_{m}$.
Since $I_{m}$ is the unique minimal ideal, we have $J=I_{m}$, therefore $%
Soc\left( G\right) =\left( 0\right) $.

Let $M$ be the unique maximal ideal of \thinspace $R$. An element $x\in R$
is invertible or zero divisor. In the last situation $x\in M$, therefore $M$
contains all zero divisors. It is clear from here that $I_{m}\subseteq M$,
since $I_{m}$ is generated by a zero divisor. Now, let $I$ be a non zero
ideal in $R$. The chain $R\supseteq I\supseteq ....\supseteq I^{\prime
}\supseteq \left( 0\right) $ is stationary, that means $I^{\prime }$ is the
minimal nonzero ideal of this chain and $I^{\prime }=I_{m}$, due to the
unicity of $I_{m}$. Therefore, $I_{m}$ is included in each nonzero ideal of $%
R$.

Assuming that $R$ is not a chain ring, then there are two nonzero ideals $I$
and $J$ such that are not included one in the other. Then we have the
following distinct chains: $\left( 0\right) \subseteq I_{m}\subseteq
...\subseteq I\subseteq R$ and $\left( 0\right) \subseteq I_{m}\subseteq
...\subseteq J\subseteq R$.  We can consider that in these chains between $R$
and $I_{m}$,  $I$ and $J$ are the last ideals strictly including $I_{m}$. If
not, we consider the last ideals strictly including $I_{m}$ from both chains
to be selected, due to Artinian ring definition.  Since, from above,
\thinspace $I\cap J\neq \left( 0\right) $ and $I\cap J$ is the minimal
nonzero ideal included in $I$ and $J$, it results that $I\cap J=I_{m}$. We
obtain that $\frac{I}{I_{m}}$ $\cap \frac{J}{I_{m}}=\widehat{I}\cap \widehat{%
J}=\left( 0\right) $ in $G$, therefore there are in $G$ two ideals $\widehat{%
I}$ and $\widehat{J}$ such that $\widehat{I}\cap \widehat{J}=\left( 0\right)
,\widehat{I},\widehat{J}\neq \left( 0\right) $, since strictly includes $%
I_{m}$. From here, we have that $\widehat{I}$ and $\widehat{J}$ are minimal
ideals in $G$. We have that $\left( 0\right) \subseteq \widehat{I}\subseteq 
\widehat{I}\oplus \widehat{J}$ and $\left( 0\right) \subseteq \widehat{J}%
\subseteq \widehat{I}\oplus \widehat{J}$ ($\widehat{I}\oplus \widehat{J}$ is
a direct sum of two proper ideals, since they are disjoint). From here,
since $\widehat{I}$ and $\widehat{J}$ are minimal ideals in $G$, we obtain $%
Soc\left( G\right) \neq \left( 0\right) $, contradiction with the fact that $%
Soc\left( G\right) =\left( 0\right) $. Therefore, we have $I\subseteq J$ or $%
J\subseteq I$ and $R$ is a chain.$\Box \medskip $

We know that a finite BL-algebra $B$ is a finite direct product of
BL-comets, $B=B_{1}\times ...\times B_{q},B_{i}$ is BL-comet. Supposing that
there is a commutative and unitary ring $R$ such that $Id\left( R\right) ~$%
has a finite BL-algebra structure, that means $Id\left( R\right)
=B_{1}\times ...\times B_{q}$. Since $Id\left( R\right) $ is finite, then $R$
is an Artinian ring and it is a finite product of Artinian local rings, $%
R=R_{1}\times ...\times R_{t}$, with $q\neq t$, then we have the following
equalities $Id\left( R\right) =Id(R_{1})\times ...\times Id(R_{t})$ and 
\begin{equation}
Id(R_{1})\times ...\times Id(R_{t})=B_{1}\times ...\times B_{q}.  \tag{2}
\end{equation}

From Proposition 25, Theorem 31 and relation $\left( 2\right) $, we can't
have $Id(R_{i})=B_{j}$, but we can have 
\begin{equation}
Id\left( R^{\prime }\right) =Id(R_{i_{1}})\times ...\times
Id(R_{i_{k}})=B_{j_{1}}\times ...\times B_{j_{s}},k\leq t,s\leq t.  \tag{3}
\end{equation}

We must remark that if $M_{i}$ is maximal ideal in $R_{i}$, then a maximal
ideal in $R$ is of the form $\mathfrak{M}_{i}=\left(
R_{1},...,M_{i},....R_{t}\right) $. The number of maximal ideals in $R$ is $t
$.  If $m_{i}$ is a minimal ideal in $R_{i}$, then a minimal ideal in $R$ is
of the form $\mathfrak{m}_{i}=\left( 0,0,...,m_{i},0,...,0\right) $. Since
each $R_{i}$ has at least a minimal ideal, the number of minimal ideals is
minimum equal with $t$. 

If all $R_{i}$ are chain rings, then $Id(R)$ is a direct product of chain
local Artinian rings, then $Id(R)$ is an MV-algebra. Therefore, in relation $%
\left( 3\right) $, we assume that at least one ring $R_{i}$ is not a chain
ring.

\textbf{Case 1.} In relation $\left( 3\right) $, we assume that at least one 
$R_{i_{j}}$ is not a chain ring, that means it has at least two minimal
ideals and one maximal ideal, from the above Lemma. It results that $%
R^{\prime }$ has at least $2k$ minimal ideals and $k$ maximal ideals. For $%
B_{j_{1}}\times ...\times B_{j_{s}}$ we have $s$ minimal ideals and at least 
$s$ maximal ideals, if all $B_{j_{i}}$ are BL-chains.

If $k<s$, then it is a contradiction with the number of maximal elements;

If $k>s$, then it is a contradiction with the numbar of minimal elements;

If $k=s$, a contradiction with the number of minimal elements.

\textbf{Case 2.} In relation $\left( 3\right) $, we assume that all $%
R_{i_{j}}$ are not chain rings, that means each of them has minimum two
minimal ideals. Then $R^{\prime }$ has at least $2k$ minimal ideals
(actually, at least  $2^{k}$) and $k$ maximal ideals. For $B_{j_{1}}\times
...\times B_{j_{s}}$ we has $s$ minimal ideals and at least $s$ maximal
ideals, if all $B_{j_{i}}$ are BL-chains.

If $k<s$, then it is a contradiction with the number of maximal elements;

If $k>s$, then it is a contradiction with the numbar of minimal elements;

If $k=s$, a contradiction with the number of minimal elements.

From the above, we obtain a contradiction and such a coomutative and unitary
ring does not exist. $\Box \medskip $

\textbf{Remark 42.} From the above Theorem, the only posibility is that $R$
to be a direct product of local Artinian rings, to each one correspond an
MV-chain, then we obtain a product of MV-chains, therefore an unordered
MV-algebra.\medskip

\textbf{Corollary 43.} \ \textit{A finite BL-algebra is a BL-comet or an
unordered MV-algebra, that means an MV-algebra which is not an MV-chain (is
a finite direct sum of MV-chains)}.

\textbf{\ \ }%
\begin{equation*}
\end{equation*}

\textbf{Conclusions.} In this paper, we studied some properties of finite
BL-comets, we gave an application of MV-algebras in cryptography, we proved
that there are no commutative and unitary rings $R$ such that its lattice of
ideals $Id\left( R\right) $ is a finite BL-algebra, which is not an
MV-algebra (Theorem 41) and we present a method to generate all BL-comets.
As a consequence, we gave a characterisation of a finite BL-algebra: it is a
BL-comet or an unordered MV-algebra. This paper closes a problem for the
study of finite BL-algebras, regarding their representation as a lattice of
ideals of commutative and unitary ring, but open a direction to study and
characterize infinte BL-algebras. \ Now, as a short notification for
readers, we must remark that even if we gave a general result in Section 3
(see Theorem 31), we also inserted a particular result (see Theorem 29) to
emphasize the way in which these results appeared. Our approach was to
consider first BL-comets of prime order, thinking at the role of the prime
numbers in the factorisation of a positive integer or in decomposition of a
finite abelian group. After that, we obtained the general result, but we
considered a good ideea to keep and present both.

\ 
\begin{equation*}
\end{equation*}

\bigskip The authors declare that there are no conflict of interests.

\begin{equation*}
\end{equation*}

\bigskip \textbf{References}%
\begin{equation*}
\end{equation*}

\textbf{[AF; 92]} Anderson, F. W., Fuller, K., (1992), \textit{Rings and
categories of modules}, Graduate Texts in Mathematics, 13(1992), 2 ed.,
Springer-Verlag, New York.

\textbf{[AM; 69]} Atiyah, M. F., MacDonald, I. G., I\textit{ntroduction to
Commutative Algebra}, Addison-Wesley Publishing Company, London, 1969.

[\textbf{BN; 09}] Belluce, L.P., Di Nola, A., \textit{Commutative rings
whose ideals form an MV-algebra}, Math. Log. Quart., 55 (5) (2009), 468-486.

[\textbf{BV;10}] Belohlavek, R., Vychodil, V., \textit{Residuated lattices
of size }$n\leq 12,$ Order, 27 (2010), 147-161.

[\textbf{CFP; 23}] C\u{a}lin, M. F., Flaut, C., Piciu, D., \textit{Remarks
regarding some Algebras of Logic}, Journal of Intelligent \& Fuzzy Systems,
45(5)(2023), Journal of Intelligent \& Fuzzy Systems, 45(5)(2023),
8613-8622, DOI: 10.3233/JIFS-232815

[\textbf{CHA; 58}] Chang, C.C.,\textit{\ Algebraic analysis of many-valued
logic}, Trans. Amer. Math. Soc. 88(1958), 467-490.

[\textbf{CFDP; 22}] Flaut, C., Piciu, D., \textit{Connections between
commutative rings and some algebras of logic}, Iranian Journal of Fuzzy
Systems, 19(6)(2022), pp. 93-110, WOS:000885481900007, DOI:
10.22111/IJFS.2022.7213,

[\textbf{FP; 22}] Flaut, C., Piciu, D., \textit{Some Examples of BL-Algebras
Using Commutative Rings}, Mathematics, 10(24)(2022), 4739, DOI:
10.3390/math10244739

[\textbf{FK; 12}] Filipowicz, M., Kepczyk, M., \textit{A note on
zero-divisors of commutative rings}, Arab J Math, 1(2012), 191--194.

[\textbf{I; 09}] Iorgulescu, A., \textit{Algebras of Logic as BCK Algebras},
A.S.E.: Bucharest, Romania, 2009.

[\textbf{NL; 03}] Di Nola, A., Lettieri, A., \textit{Finite BL-algebras},
Discrete Mathematics 269(2003), 93 -- 112.

[\textbf{NL; 05}] Di Nola, A., Lettieri, A., \textit{Finiteness based
results in BL-algebras}, Soft Comput 9(2005) 9, 889--896, DOI
10.1007/s00500-004-0447-7.

[\textbf{P; 07}] Piciu, D., \textit{Algebras of fuzzy logic}, Ed.
Universitaria, Craiova, 2007.

[\textbf{TT; 22}] Tchoffo Foka, S. V., Tonga, M., \textit{Rings and
residuated lattices whose fuzzy ideals form a Boolean algebra}, Soft
Computing, 26 (2022) 535-539.

[\textbf{TT; 17}] Thakur, K., Tripathi, B.P., \textit{A Variant of NTRU with
split quaternions algebra}, Palestine Journal of Mathematics, 6(2)(2017),
598--610.

[\textbf{T; 99}] Turunen, E., \textit{Mathematics Behind Fuzzy Logic},
Physica-Verlag, 1999.

[\textbf{WD; 39}] Ward, M., Dilworth, R.P., \textit{Residuated lattices},
Trans. Am. Math. Soc. 45(1939), 335--354.%
\begin{equation*}
\end{equation*}

Cristina Flaut$^{\ast }$(corresponding author)

{\small Faculty of Mathematics and Computer Science, Ovidius University,}

{\small Bd. Mamaia 124, 900527, Constan\c{t}a, Rom\^{a}nia,}

{\small \ http://www.univ-ovidius.ro/math/; https://www.cristinaflaut.com,}

{\small e-mail: cflaut@univ-ovidius.ro; cristina\_flaut@yahoo.com}

\bigskip

Dana Piciu

{\small Faculty of \ Science, University of Craiova, }

{\small A.I. Cuza Street, 13, 200585, Craiova, Rom\^{a}nia,}

{\small http://www.math.ucv.ro/dep\_mate/}

{\small e-mail: dana.piciu@edu.ucv.ro, piciudanamarina@yahoo.com}

\bigskip

Bianca Liana Bercea-Straton

{\small PhD student at Doctoral School of Mathematics,}

{\small Ovidius University of Constan\c{t}a, Rom\^{a}nia}

{\small e-mail: biancaliana99@yahoo.com}

\end{document}